\documentclass[final,onefignum,onetabnum]{siamart250211}

\usepackage{amssymb}
\usepackage{mathrsfs} 
\usepackage{lipsum}
\usepackage{amsfonts}
\usepackage{mathtools}
\usepackage{graphicx}
\usepackage{subcaption}
\usepackage{booktabs} 
\usepackage{epstopdf}
\usepackage{tikz}
\usepackage{amsmath}
\usepackage{algorithmic}
\usepackage{algorithm}     
\usepackage{extarrows}
\ifpdf
  \DeclareGraphicsExtensions{.eps,.pdf,.png,.jpg}
\else
  \DeclareGraphicsExtensions{.eps}
\fi

\newsiamremark{remark}{Remark}
\newsiamremark{hypothesis}{Hypothesis}
\crefname{hypothesis}{Hypothesis}{Hypotheses}
\newsiamthm{claim}{Claim}
\newsiamthm{assumption}{Assumption}      
\crefname{assumption}{Assumption}{Assumptions}  
\newsiamthm{example}{Example}
\crefname{example}{Example}{Example}

\headers{An inverse random diffraction grating problem}{Zhiqi Sun  and Yiwen lin}

\title{An inverse random diffraction grating problem for the Helmholtz equation\thanks{Submitted to the editors Jul 22, 2025.
\funding{This work was partially supported by National Key R\&D Program of China (no. 2024YFA1016100) and National Natural Science Foundation of China (no. U21A20425, no. 12201404). }}}

\author{Zhiqi Sun\thanks{School of Mathematical Sciences, Shanghai Jiao Tong University, Shanghai 200240, China
 (\email{sunzhq1016@sjtu.edu.cn}).}
\and 
Yiwen lin\thanks{Corresponding author. School of Mathematical Sciences, Shanghai Jiao Tong University, Shanghai 200240, China
  (\email{linyiwen@sjtu.edu.cn}).}}

\usepackage{amsopn}

\makeatletter
\newcommand*{\addFileDependency}[1]{
  \typeout{(#1)}
  \@addtofilelist{#1}
  \IfFileExists{#1}{}{\typeout{No file #1.}}
}
\makeatother

\newcommand*{\myexternaldocument}[1]{%
    \externaldocument{#1}%
    \addFileDependency{#1.tex}%
    \addFileDependency{#1.aux}%
}

\usepackage{soul}  
\setstcolor{magenta}

\ifpdf
\hypersetup{
  pdftitle={Stochastic Grating Profiles},
  pdfauthor={D. Doe, P. T. Frank, and J. E. Smith}
}
\fi

\myexternaldocument{ex_supplement}

\begin{document}

\maketitle

\begin{abstract}This paper investigates the inverse scattering problem of time-harmonic plane waves incident on a perfectly reflecting random periodic structure. To simulate random perturbations arising from manufacturing defects and surface wear in real-world grating profiles, we propose a stochastic surface modeling framework motivated by the discretization of the Wiener process. Our approach introduces randomness at discrete nodes and then applies linear interpolation to construct the surface, marking a novel attempt to incorporate the concepts of the Wiener process into random surface representation. Under this framework, each realization of the random surface generates a Lipschitz-continuous diffraction grating, mathematically represented as a sum of a baseline profile and a weighted linear combination of local `tent' basis functions, meanwhile preserving key statistics of the random surface.  Building on this representation, we introduce the Recursive Parametric Smoothing Strategy (RPSS) to invert the key statistics of our random surfaces. Combined with Monte Carlo sampling and a wavenumber continuation strategy, our reconstruction scheme demonstrates effectiveness across multiple benchmark scenarios. Several numerical results are presented along with some discussions in the end on reconstruction mechanisms and future extensions.

\end{abstract}



\section{Introduction}

The inverse diffraction grating problem is a classic and challenging topic at the intersection of several major branches of mathematics. With wide-ranging applications in fields such as geophysical prospecting and medical imaging, inverse scattering holds a crucial position not only in mathematical scattering theory but also in modern diffractive optics. Consider a two-dimensional periodic perfectly conducting gratings, where a time-harmonic electromagnetic plane wave incidents on this structure and subsequently scatters above this surface. Assuming the medium is invariant in the $z$-direction, the  three-dimensional time-harmonic Maxwell equations will be reduced to a two-dimensional Helmholtz equation. For a comprehensive overview of scattering problems in periodic structures, we refer the reader to \cite{bao2022maxwell} and the references therein.

For wave scattering by periodic structures, the direct (or forward) problem is to determine the scattered field from a known incident wave and surface profile. Conversely, the inverse problem entails reconstructing the diffractive surface from measured scattered patterns. Extensive research has addressed both forward and inverse problems, spanning theoretical investigations  to computational methodologies. For the forward problems where the surface is  $\mathcal{C}^2$-smooth,  theoretical uniqueness results of the scattering field can be found in \cite{kirsch2005diffraction}. Computationally, prominent  numerical methods primarily include the integral equation approach \cite{nedelec1991integral}, variational approach \cite{bruno1993numerical} and adaptive finite-element method \cite{bao2005adaptive}.
For the inverse problems, \cite{hettlich1997schiffer,kirsch1994uniqueness,bao1994uniqueness,bao1995inverse} established the uniqueness of reconstructing smooth surface from scattering data within their corresponding respective problem settings. In parallel, practical and effective methods have been developed, such as high-order perturbation approach \cite{ito1999high}, and factorization method \cite{arens2003factorization}. Relaxing these strict smoothness conditions to include Lipschitz continuous gratings has also received attention in the literature. For forward problems, \cite{zhang2004integral} established the uniqueness of scattering solutions represented by a combination of single- and double-layer potentials, while \cite{elschner2002inverse} achieved similar results using only single-layer representations. The corresponding inverse problem, however, presents a greater challenge; the global uniqueness question for general Lipschitz surfaces remains largely unresolved. A notable contribution was made by \cite{bao2011unique}, where uniqueness was proven for the specific case of two-dimensional periodic polyhedral Lipschitz structures under single-frequency incidence. Consequently, most existing studies have focused primarily on developing computational approaches for reconstructions \cite{bruckner2003two, bao2012computational, zhang2015efficient} . Further relevant results will be discussed where appropriate in subsequent sections.

 Nowadays, the widespread existence of uncertainties has spurred growing research interest in inverse scattering problems. Since practical grating structures inherently possess manufacturing defects, surface wear and operational damage during daily use, all of which contribute to a considerable degree of uncertainty, incorporating randomness in surface scattering problems will more closely align with real-world problem scenarios. However, this field remains largely limited to investigations of random source problems like \cite{bao2014inverse,li2011inverse}, with other aspects of uncertainty receiving comparatively little attention due to the substantial mathematical and computational challenges posed by random surface geometries and their often inexplicit functional forms. In this work, we conduct a preliminary exploration of the inverse random surface scattering problem: Given measured scattering wave fields generated by time-harmonic plane wave incidence, we aim to determine key statistical properties of the periodic random structure from data acquired  at certain $y=y_0$ above the surface. To the best of our knowledge, the study of random surfaces is still in its infancy, with only initial progress having been made. Monte Carlo method was used to compute the expected scattering patterns of random rough surfaces in  \cite{johnson1996backscattering}. An ultrasonic methodology was adopted in \cite{shi2017recovery} to reconstruct the height correlation function of remotely inaccessible random rough surfaces in solids . More recently, Bao, Lin and Xu\cite{bao2020inverse} proposed the MCCUQ algorithm, which integrates the Landweber iteration method with Uncertainty Quantification (UQ) techniques. This innovative approach has shown promising experimental results in reconstructing the statistical parameters of random surfaces, highlighting its potential for practical applications. Subsequently, random periodic structures were studied in other contexts,  such as for the elastic wave equation in \cite{gu2024inverse}. The authors combined the modified Monte Carlo continuation method and kernel density estimation method to reconstruct key statistical properties,  demonstrating effectiveness for both Gaussian and non-Gaussian random surfaces. Building on these approaches, Lv, Wang \itshape{et al.} \normalfont studied the random penetrable periodic structures in \cite{wang2025numerical} and the acoustic-elastic interaction with random periodic interfaces in \cite{wang2025numericall}. Beyond these specialized studies, research on inverse problems involving random surfaces remains remarkably scarce, particularly for cases without strict smoothness assumptions.

In this paper, we investigate a specific class of random periodic surfaces with relaxed smoothness constraints. Specifically, we adopt a formal stochastic representation for surfaces inspired by the random source modeling in  \cite{bao2014inverse,li2011inverse}:
\begin{equation}{\label{form}}
    ``f(\omega,x)=g(x)+h(x)dW_x(\omega)".
\end{equation}

Here $g(x)$, $h(x)$ are deterministic real functions with compact support, $W_x$ is a one-dimensional spatial Wiener process and $dW_x$ is its stochastic differential in the It$\hat{\text{o}}$ sense, commonly used to model white noise. In the context of our surface scattering problem, this formulation admits a clear physical interpretation: the function   $g(x)$ captures the  deterministic, macroscopic baseline profile, representing the macroscopic surface geometry; the function  $h(x)$ controls the  perturbation amplitude by spatially modulating the noise intensity while the  ``$dW$'' term encodes the surface roughness. However, since  white noise  $dW_x$  is a generalized function (or distribution), the formal expression \cref{form}  lacks rigorous meaning as a pointwise-defined univariate function in the strict mathematical sense.  To address this, Bao, Li \it{et al.} \normalfont rigorously reformulated   \cref{form} as the stochastic integral $ f(x)=g(x)+\int_0^xh(s)dW_s(\omega)$ which transforms the white noise into a continuous stochastic process suitable for numerical approximation. Motivated by their integral-based representation and discrete nodal sampling strategy, we adopt the same discrete nodal sampling strategy to preserve nodal variance consistency. Instead of using cumulative noise increments which would introduce inter-node correlations and analytical complexities, our model utilizes independent Gaussian sampling at each node, followed by linear interpolation to construct Lipschitz-continuous  surfaces, which not only ensure subsequent numerical analysis, but also preserve key statistical properties of the scattering surface. This finite-dimensional, locally independent model can be viewed as a Galerkin projection of  $dW_x$ onto the subspace of piecewise linear basis functions.  It  aims to simulate scenarios involving locally independent random disturbances, such as localized vibrations or surface defects. To enhance physical fidelity, this paradigm will be further refined by two key assumptions detailed later. Building upon these foundational considerations, we develop the Recursive Parametric Smoothing Strategy (RPSS) as a novel inversion framework for characterizing the statistics of random grating profiles. By integrating RPSS with Monte Carlo sampling and a wavenumber continuation approach, we propose the Monte Carlo Continuation $h(x)$  (MCCh) algorithm, realizing effective reconstructions across diverse benchmark configurations, particularly when $h(x)$ is expressed by the linear combination of elementary trigonometric functions. In contrast with \cite{bao2020inverse} which modeled the random surface by the Karhunen-Lo\`eve (KL) expansion and focused on analyzing coefficients of known KL basis functions, our work models the random surface motivated by the Wiener process with relaxed smoothness constraints and shifts the emphasis to inverting the functional forms of the basis modes themselves. 
Thus, this work can also be regarded as an extension of the study in \cite{bao2020inverse}.

The remainder of this paper is organized as follows. We begin in \Cref{sec:For} by establishing the physical modeling framework and mathematical problem formulation. In \Cref{inv}, we present the proposed inversion methodology and algorithmic implementation. We then validate the method through a series of numerical benchmark cases in  \Cref{sec:experiments}. \Cref{sec:obs} provides an analysis towards experiment results, along with other related discussions. Finally, \Cref{sec:conclusions} concludes the paper by summarizing our findings, discussing inherent limitations, and outlining directions for future research.

\section{Problem Formulation}
\label{sec:For}
We shall start by presenting our modeling approach for random surfaces in \Cref{sub1} and subsequently delve into the forward and inverse modeling problems in \Cref{sub2}.
\subsection{Physical modeling of random surfaces}\label{sub1}

\begin{figure}[htbp]
    \centering
    \begin{subfigure}{0.33\linewidth}
        \includegraphics[width=\linewidth]{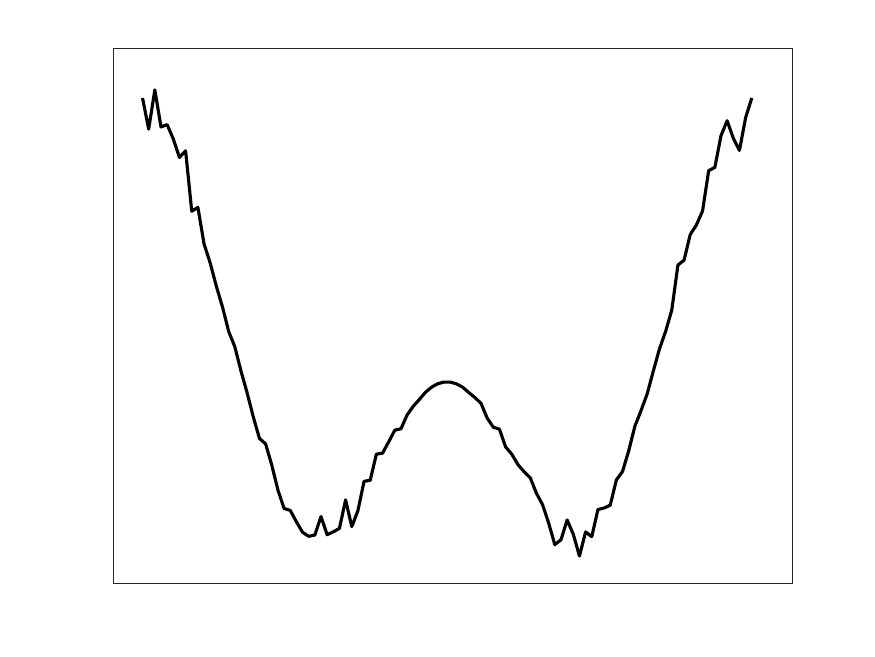}
        \caption{Realization 1}
    \end{subfigure}
    \begin{subfigure}{0.33\textwidth}
        \includegraphics[width=\linewidth]{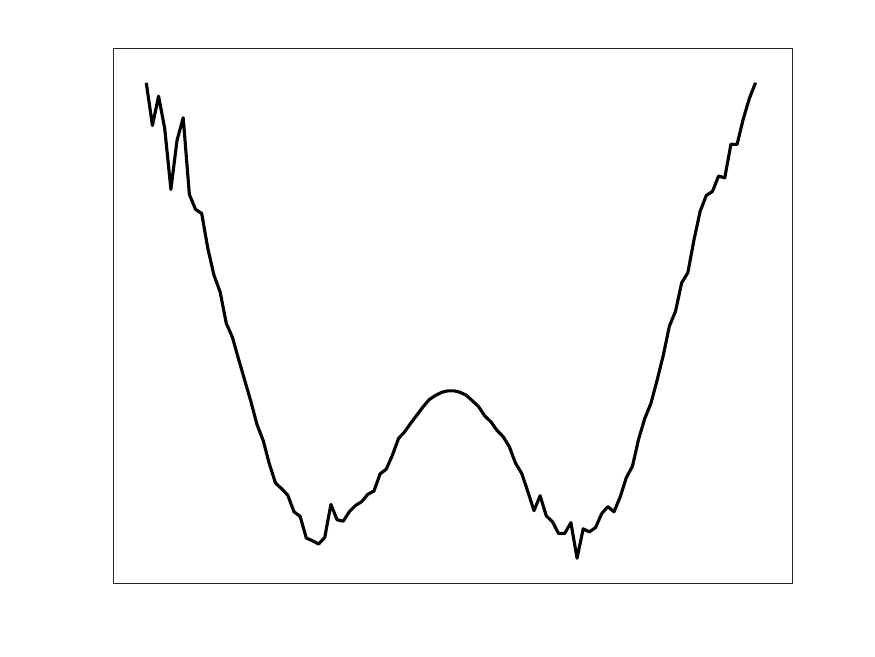}
        \caption{Realization 2}
    \end{subfigure}
    
    
    \begin{subfigure}{0.33\textwidth}
        \includegraphics[width=\linewidth]{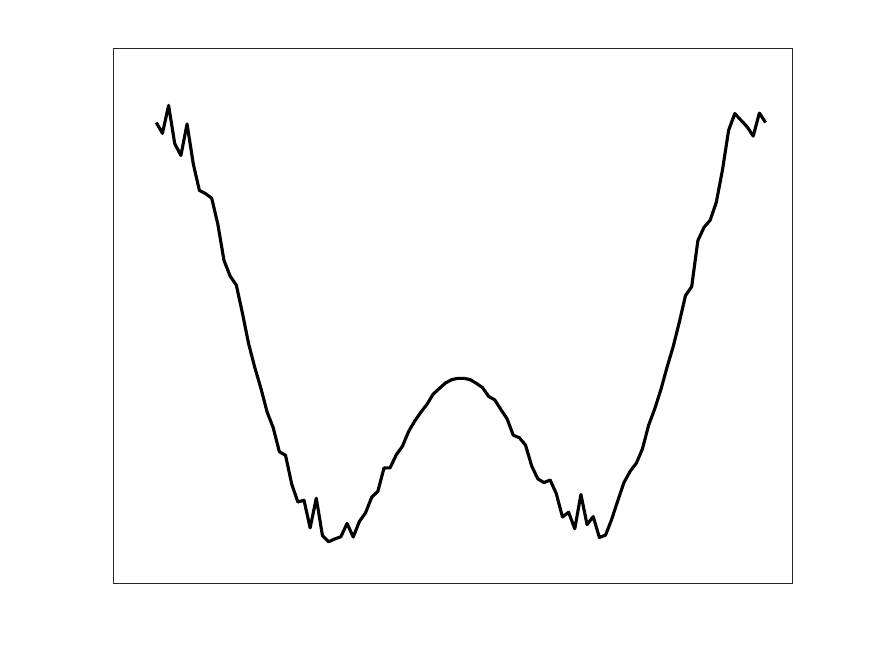}
        \caption{Realization 3}
    \end{subfigure}
    \begin{subfigure}{0.33\textwidth}
        \includegraphics[width=\linewidth]{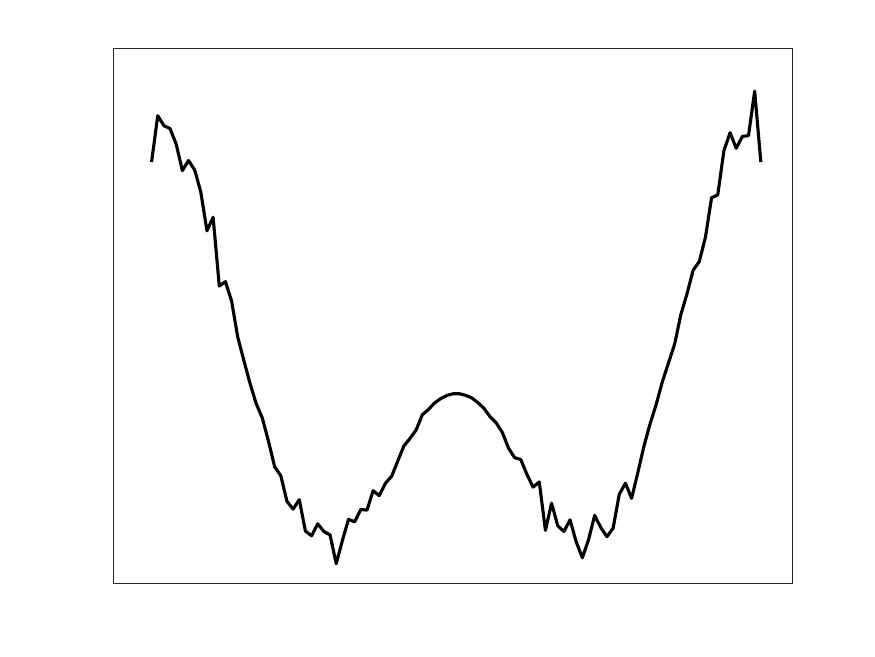}
        \caption{Realization 4}
    \end{subfigure}
    
    \caption{Four realizations of one example of the random surface \cref{ourform}}
    \label{fig:four_images}
\end{figure}

 The profile of diffraction grating is represented by
\[\Gamma_f:=\{ (x,y)\in\mathbb{R}^2|y=f(\omega,x)\},\]
defined in a probability space  $(\Omega,F,\mu)$ where  $\omega\in\Omega$ denotes the random sample, $(x,y)\in\mathbb{R}^2$. Each realization of the random surface  $f_m(x)=f(\omega_m,x)$  can be viewed as a one-dimensional real function over the spatial interval $[0,\Lambda]$, where $\Lambda$ denotes one period.  Our modeling philosophy aligns with  the inverse random source literatures \cite{bao2014inverse,li2011inverse} where the random source terms are represented as the formal expression $f(\omega,x) = g(x) + h(x)dW_x(\omega)$. While physically intuitive, this notation  lacks classical mathematical rigor as a pointwise-defined function due to the generalized nature of the spatial white noise $dW_x$, which exists strictly as a generalized random field or distribution \cite{sobczyk2013stochastic}. To address this rigorously, Bao, Li \it{et al.} \normalfont developed a mathematically framework by interpreting the stochastic term $h(x)dW_x$ as an 
It$\hat{\text{o}}$ integral with respect to spatial Brownian motion $\{W_x\}_{x\geqslant0}$:
\begin{equation}\label{pao}
   f(\omega,x) = g(x) + h(x)dW_x(\omega)\longrightarrow f(\omega,x)=g(x)+\int_0^xh(s)dW_s(\omega).
\end{equation}
The integral in \cref{pao} converts the highly irregular white noise into a continuous stochastic process, thereby producing a well-defined classical one-dimensional function suitable for numerical approximation and theoretical analysis. Following this integral-based representation, they adopt discrete node sampling
\begin{equation}{\label{paoo}}
    \int_{0}^{1} h(y) \, dW_{y} \approx \sum_{j=0}^{N_0-1} h(x_j) \Delta W_j,\quad \text{where } \Delta W_j = \xi_j \sqrt{\Delta x}, \  \xi_j \overset{\text{i.i.d.}}{\sim} \mathcal{N}(0,1)
\end{equation}
to approximate the integral numerically, which inspires our modeling approach of the random surface.  To reduce computational complexity while maintaining accuracy and enabling direct numerical analysis, we propose an alternative finite-dimensional model for random surfaces that explicitly captures independent perturbations emerging on surfaces at individual spatial locations. This approach is particularly suited for modeling localized surface vibrations or isolated disturbances that behave as independent random sources. Specifically, our construction involves sampling independent Gaussian random variables $\xi_i$ at discrete spatial nodes $x_i\in[0,\Lambda]$, then modulated by the intensity function $h(x_i)$ as:
\begin{equation}{\label{paooo}}
    f(x_i)=g(x_i)+h(x_i)\xi_i\sqrt{\Delta x},\quad i\in\{1,2,\cdots,N_0\},
\end{equation}
and subsequently interpolating these node values linearly. Our node-based sampling strategy aligns closely with the sampling approach used in \cite{li2011inverse} and \cite{bao2014inverse} where values at discrete points are similarly generated to approximate spatial integrals of white noise. Here the nodal points $\{x_i\}_{i=1}^{N_0}$ are determined by the hyperparameter $\Delta x$ since $\Delta x=x_i-x_{i-1}=\frac{\Lambda}{N_0}$, where $N_0$ denotes the number of sample points spanning the interval $[0,\Lambda]$. This construction finally yields Lipschitz-continuous paths which maintains the statistical characteristics of the unknown scattering profile and offers significant numerical and analytical advantages. The complete surface modeling procedure comprises four key steps:



\begin{enumerate}
  \item \textbf{Spatial Discretization:} The interval $[0,\Lambda]$ is partitioned into $N_0$ equal segments, with each segment having a length of $\Delta x$
    \[0=x_0<x_1<\cdots<x_N=\Lambda,\quad \Delta x=x_i-x_{i-1}=\frac{\Lambda}{N_0}.\]
    Note that $N_0$ or $\Delta x$ is an essential hyper-parameter in our setting;
    \item \textbf{Modeling of Wiener increments:} Approximate the continuous Wiener process $W(x)$ by independent normal increments on each subinterval
    \[dW_i=W(x_i+\Delta x)-W(x_{i})\sim\mathcal{N}(0,\Delta x).\] 
    In this way, we obtain a stochastic sequence $\{\Delta W_i\}$ defined on the grid points;
    \item \textbf{Approximating $dW$ at grid points:} For each surface, we assume $W(x_0)=0$. Then at the $i$-th grid point $x 
_{i}$, $i\geqslant1$, we define:
\[f(\omega,x_i)=g(x_i)+h(x_i)\xi_{i}\sqrt{\Delta x},\]
where $\{\xi_{i}(\omega)\}_{i=1,\cdots,N_0}\overset{\text{i.i.d}}{\sim}\mathcal{N}(0,1)$ are random variables;
\item \textbf{Interpolation-based surface construction:}
Apply linear interpolation between adjacent nodes $x_i$ and $x_{i+1}$
\[ f(\omega, x) = \frac{x_{i+1} - x}{\Delta x} f(\omega, x_i) + \frac{x - x_i}{\Delta x} f(\omega, x_{i+1}), \quad x_i < x < x_{i+1}, \]
then we are able to  obtain a Lipschitz continuous curve $f(\omega, x)$  reflecting white noise perturbations while being numerically tractable.
\end{enumerate}
Therefore, our complete random surface model takes the final form:
\begin{equation}{\label{ourform}}
  f(\omega,x) = g(x) + \sum_{j=0}^{N_0} h(x_j) \xi_j(\omega) \sqrt{\Delta x} \varphi_j(x),
\end{equation}
where $\varphi_j(x)$ is the piece-wise linear ``tent'' function centered at $x_j$:
\[
\varphi_j(x) = 
\begin{cases} 
1 - \dfrac{|x - x_j|}{\Delta x}, & |x - x_j| \leq \Delta x, \\
0, & |x - x_j| > \Delta x.
\end{cases}
\]
An example of our random surface is shown in \cref{fig:four_images}.
The stochastic term in our model \cref{ourform} can be interpreted as a finite-dimensional Galerkin projection of the white noise process $dW_x$ onto the subspace spanned by piecewise linear basis functions. In contrast to globally correlated representations like the Karhunen-Lo\`eve (KL) expansion,  our model considered in this paper emphasizes local independence. A key feature of our finite-dimensional representation is the hyperparameter $\Delta x$, which  controls both the spatial resolution of our random surface and the variance magnitude at each sampled node. As $\Delta x\rightarrow 0$,the variance diminishes, and the Lipschitz continuous random surface $f^m_{lip}(x)$ degenerates towards the deterministic function $g(x)$. Besides this, an excessively small $\Delta x$ may introduce overly frequent and dense perturbations, leading to high-frequency oscillations, significant surface distortion and numerical instabilities. Therefore, to ensure practical applicability, numerical stability and computational efficiency, we  propose the following assumption:

\begin{assumption}[Lower Bound on Spatial Discretization]{\label{110}}There exists a small positive lower bound $\varepsilon>0$  such that the   interval length $\Delta x$ is no smaller than $\varepsilon$.
\end{assumption}
This assumption guarantees that the random perturbations have moderate spatial oscillations, aligning the modeled surface more closely with realistic physical scenarios. In our   numerical presentations we will investigate and discuss how varying $\Delta x$ affects computational stability and accuracy, aiming to balance model simplicity with reconstruction precision. For the simulations presented in this paper, we set this lower bound to  $\varepsilon\approx0.004$ (i.e., $N_0$=160).

Furthermore, we  assume that the intensity of the random perturbation remains moderate relative to the deterministic surface profile $g(x)$. Excessively large perturbations could obscure the essential physical characteristics encoded in $g(x)$, thereby compromising meaningful analysis. To formalize this constraint, we introduce the following assumption:

\begin{assumption}[Moderate Amplitude of Random Perturbation]{\label{911}}
There exists a constant $c\in (0,1),$ such that for any random realization $f(\omega,x)$, the following inequality holds:
\[\Big\|f(\omega,x)-g(x)\Big\|_{\mathcal{L}^2}\leqslant c\Big\|g(x)\Big\|_{\mathcal{L}^2},\]
where $\|\cdot\|_{\mathcal{L}^2}$ represents the $\mathcal{L}^2$ norm.
\end{assumption}

With \cref{110} and \cref{911} established, we now proceed to the main problem formulation.

\subsection{Direct and inverse diffraction problems}\label{sub2}
In this part we shall introduce the problem formulation. Suppose a  plane wave incident on a periodic random surface $\Gamma_f$. With all the preparations in \Cref{sub1},  each realization of the random surface $f(\omega_m,x):=f^{m}_{lip}(x)\in \mathcal{C}_p^{0,1}$, i.e., $f^{m}_{lip}(x)$ is a periodic Lipschitz function of period $\Lambda=2\pi$. The space below gratings is filled with some perfectly reflecting material. Let $\Omega_f:=\{(x,y)\in\mathbb{R}^2:y>f(x),x\in\mathbb{R}\}$ be filled with homogeneous material whose wavenumber $\kappa>0$. The incident plane wave has the form $u^i=e^{i\alpha x-i\beta y}$, $\alpha=\kappa sin\theta$, $\beta=\kappa cos\theta$, 
and $\theta\in(-\frac{\pi}{2},\frac{\pi}{2})$ is the incident angle with respect to $y$-axis. Define the total field $u^{\text{total}}=u^i+u^s(\omega;\cdot)$ in $\Omega_f$, with $u^s$ representing the scattered field. Then the scattered field in the TE (transverse electric) mode satisfies the Helmholtz equation 
\begin{equation}\label{23}
    \Delta u^s(\omega;\cdot)+\kappa^2u^s(\omega;\cdot)=0,\quad\text{in}\ \Omega_f,
\end{equation}
with a Dirichlet boundary condition:
\begin{equation}\label{24}
    u^s(\omega;\cdot)+u^i=0,\quad\text{on}\ \Gamma_f.
\end{equation}
Due to uniqueness consideration, we seek for the quasi-periodic solutions which satisfy
\begin{equation}\label{25}
    u^s(\omega; x+\Lambda, y) = e^{i\alpha\Lambda} u^s(\omega; x, y), \quad\text{ in } \Omega_f.
\end{equation}
Moreover, we require $u^s$ to satisfy a radiation condition, i.e., $u^s$ admits a Rayleigh expansion:\begin{equation}\label{26}
    u^s(\omega; \cdot) = \sum_{n \in \mathbb{Z}} A_n e^{i \alpha_n x + i \beta_n y}, \quad y > \max_{x \in (0, \Lambda)} f(\omega; x),
\end{equation}
where
\begin{equation}\label{27}
    \alpha_n := \alpha + \frac{2\pi n}{\Lambda}, \quad
\beta_n := \left\{
\begin{array}{ll}
\sqrt{\kappa^2 - \alpha_n^2}, & |\kappa| > |\alpha_n|, \\
i\sqrt{\alpha_n^2 - \kappa^2}, & |\kappa| < |\alpha_n|.
\end{array}
\right.
\end{equation}
So the direct problem is to determine scattering field $u^s$ based on the information of incidence field $u^i$ and surface $f(\omega,x)$ where $u^s$ satisfied \cref{23}-\cref{26}. For general Lipschitz grating profiles, the existence of a unique solution to the Dirichlet problem \cref{23}-\cref{27} is established in \cite{elschner2002inverse}. Conversely, the inverse problem is to reconstruct the key statistical properties of the random surface, namely the mean profile $g(x)$ and  variance intensity $h(x)$, from the incident wave $u^i$ and measurements of the scattered field $u^s$.




\section{Reconstruction Method}{\label{inv}}
For each fixed sample $\omega_m$, the surface is modeled as a deterministic Lipschitz continuous grating under our problem formulation. We aim to reconstruct $g(x)$ and $h^2(x)$ with all the sample data of incident wave field $u^i$ and the multiple frequency scattering field $u^s(\omega,x,y_0)$. The scattering data $u^s(\omega,x,y_0)$ are measured at an artificially selected boundary $\Gamma_0=\{(x,y)\in\mathbb{R}^2|y=y_0\}$ above $\Gamma_f$ where  $y_0$ is chosen such that it is greater than the maximum possible height of the random surface.

This section details our reconstruction methodology. We begin in \Cref{sec:alg} by reviewing two established inversion schemes to provide context. The problem is then decomposed into an ill-posed linear subproblem in \Cref{32} and a well-posed nonlinear subproblem in \Cref{33}. Our proposed Recursive Parametric Smoothing Strategy (RPSS) is introduced in \Cref{RPSS}. Following this, the numerical implementation using the Monte Carlo Continuation (MCC) method and its modification are detailed in \Cref{MCCC} and \Cref{MC}, respectively. Finally, \Cref{acc} presents the complete algorithm for inverting the key statistical quantities of the random diffraction grating in this work.

\subsection{Two different inversion schemes of density function and surface}
\label{sec:alg}
The inversion method requires processing the entire set of data samples, which first necessitates a method for handling each individual realization. We characterize the scattered field $u^s$ through the single-layer potential representation:
\begin{equation}\label{90}
    u^s(\omega; x, y) = \int_0^\Lambda \varphi(\omega; s) G(x, y; s, 0) \, \mathrm{d}s,
\end{equation}
where $\varphi\in\mathcal{L}^2$ is an unknown periodic density function and $G$ is  the quasi-periodic Green function:
\begin{equation}\label{91}
    G(x, y; s, t) = \frac{i}{2\pi} \sum_{n \in \mathbb{Z}} \frac{1}{\beta_n} e^{i\alpha_n (x-s) + i\beta_n |y-t|}, \quad (x, y) \neq (s, t).
\end{equation}
For a given deterministic surface, once we have the measured scattering data $u^s(x,y_0)$,  the expression of density function $\varphi(x)$ can be calculated, thereby enabling the determination of surface $f(x)$.
For Lipschitz continuous surfaces, two primary schemes exist to recover the density function $\varphi(x)$ and the surface profile $f(x)$. The first  is  the ``decomposition scheme", which splits the problem into a linear severely ill-posed subproblem and a nonlinear well-posed subproblem. In this approach, the computational workflow begins with solving  ill-posed subproblem through  regularization techniques,  yielding an explicit expression for $\varphi$. This reconstructed  $\varphi(x)$ is then incorporated into the subsequent well-posed nonlinear subproblem to iteratively solve for  $f(x)$ \cite{bruckner2003two}; In contrast, the second scheme maintains $\varphi(x)$ in an implicit form, solving for both $\varphi(x)$ and $f(x)$ simultaneously within a coupled reconstruction framework at the boundary  $\Gamma_f$ \cite{elschner2002inverse}. Given the availability of multi-frequency scattering data and the need to reduce computational costs, our work adopts the decomposition framework. 

\subsection{The linear ill-posed  problem}{\label{32}}
In this part, we derive the components needed for our inversion algorithm. Following the decomposition framework, the linear ill-posed subproblem is addressed first. Our initial goal is to establish a formal relationship between the unknown density function $\varphi(x)$ and the measured scattered field $u^s$. Taking $y=y_0$ into \cref{90} yields: $u^s(\omega; x, y_0) = \int_0^\Lambda \varphi(\omega; s) G(x, y_0; s, 0) \, \mathrm{d}s$. Expand $u^s$ as:
\begin{equation}\label{88}
    u^s(\omega; x, y_0) = \sum_{n \in \mathbb{Z}} u_n(\omega) e^{i \alpha_n x},\quad   u_n(\omega) = \frac{1}{\Lambda} \int_0^\Lambda u^s(\omega; x, y_0) e^{-i \alpha_n x} \, \mathrm{d}x,
\end{equation}
due to its quasi-periodicity. This allows us to derive an explicit formulation for  the density function $\varphi$:
\begin{equation}\label{shiit}
    \varphi(\omega; s) = \sum_{n \in \mathbb{Z}} \varphi_n(\omega) e^{i \alpha_n s},\quad \varphi_n(\omega) = -i \frac{2\pi}{\Lambda} \beta_n u_n(\omega) e^{-i \beta_n y_0}.
\end{equation}
While this equation provides a direct link between $\varphi$ and the data $u^s$, it is numerically unstable for evanescent modes  due to the exponential growth of the $\beta_n$ term which would amplify noise. To address this, we implement Tikhonov regularization scheme \cite{engl1996regularization} to stabilize the inversion. As evidenced by \cite{bao2012computational,bao2020inverse}, the Tikhonov regularization approach achieves significant computational advantages, eliminating the need for intricate solution procedures.  
   For the subsequent nonlinear problem of finding the surface $f(x)$, define operator $\mathcal{T}_f$: $\mathcal{L}^2(\Omega\ ;\mathcal{L}^2(0,\Lambda))\rightarrow \mathcal{L}^2(\Omega\ ;\mathcal{L}^2(0,\Lambda))$ by:
\begin{align}
  &\mathcal{T}_f [\varphi](\omega; x) = \int_0^\Lambda \varphi(\omega; s) G(x, f(x); s, 0) \, \mathrm{d}s.
  \label{010}
\end{align}  
which maps the density function $\varphi$ to the field on the boundary  $\Gamma_f$. 
Combining \cref{88} and \cref{shiit} together with regularization techniques, $\mathcal{T}_f$ culminates in an explicit formulation:

\begin{equation}\label{heihei}
    \mathcal{T}_f [\varphi](\omega; x) = \sum_{n \in \mathbb{Z}} \psi_n(\omega) e^{i \alpha_n x + i \beta_n f(\omega; x)},
\end{equation}
where
\begin{equation}\label{reg}
    \psi_n(\omega) = 
\begin{cases} 
u_n(\omega) e^{-i \beta_n y_0}, & \text{for } \kappa > |\alpha_n|, \\
u_n(\omega) \frac{e^{i \beta_n y_0}}{e^{2 i \beta_n y_0} + \gamma}, & \text{for } \kappa < |\alpha_n|,
\end{cases}
\end{equation}
with $\gamma$  a positive regularization parameter.

\subsection{The nonlinear well-posed  problem}{\label{33}}
With the components for the forward operator $\mathcal{T}_f$ established, we now address the nonlinear, well-posed problem of reconstructing the surface profile $f$. This nonlinear problem arises form  the Dirichlet boundary condition \cref{24} which mathematically requires the total field vanishing at the boundary $\Gamma_f$. With the notations above, \cref{24} will turn into:
\begin{equation}\label{heii}
 \Big\| (\mathcal{T}_f [\varphi])(\omega; x) + u^i(x, f(\omega,x)) \Big\|^2_{\mathcal{L}^2(0, \Lambda)} = 0,
\end{equation}
By substituting the concrete form of $\mathcal{T}_f$ computed in \cref{heihei} into \cref{heii}, one has
\begin{equation}\label{eq:01}
    \left\| \sum_{n \in \mathbb{Z}} \psi_n e^{i \alpha_n x + i \beta_n f(\omega,x)} + e^{i \alpha x - i \beta f(\omega,x)} \right\|^2_{\mathcal{L}^2(0, \Lambda)} = 0.
\end{equation}
For practical implementation, due to the exponential decay of the Rayleigh mode amplitudes $\psi_n$ with increasing $n$, we only need to take a sufficiently large truncation order $N$ into consideration, which will lead to
\begin{equation}\label{eq:02}
    \left\| \sum_{n = -N}^{n=N} \psi_n e^{i \alpha_n x + i \beta_n f(\omega,x)} + e^{i \alpha x - i \beta f(\omega,x)} \right\|^2_{\mathcal{L}^2(0, \Lambda)} = 0.
\end{equation}

\subsection{Recursive Parametric Smoothing Strategy}{\label{RPSS}}
Our ultimate goal is to reconstruct statistics of the random surface $f(\omega,x)$. This requires solving  the nonlinear equation \cref{eq:02} for each realization where $f(x)$ takes the Lipschitz-continuous form $f^m_{lip}(x)$. Although each realization of the random surface admits an explicit mathematical representation, the non-parametric and non-smooth nature of $f^m_{lip}(x)$ presents substantial computational challenges in solving \cref{eq:02}. 
To overcome this, prior work has shown the effectiveness of parameterizing grating profiles using Fourier series, such as \cite{chang2024novel}, which establishes that Fourier series representations permit effective parameter characterization and iterative updating, whereas \cite{bruckner2003two} develops alternative parametrization schemes for binary and piece-wise linear grating structures. For our problem, the primary objective is to recover the statistical functions  $g(x)$ and $h^2(x)$ rather than exact surface topology of $f_{lip}^m(x)$. Therefore, we require a parameterization that is both computationally efficient and sufficiently expressive to capture these statistics. These considerations motivate our  adoption of Fourier-series parametrization, where we approximate the true surface $f_{lip}^m(x)$ with a smooth, infinitely differentiable function:
\begin{equation}{\label{smoo}}
   f(\omega_m,x):= f_{m,\infty}(x)=c_{m,0}+ \sum_{p\in\mathbb{N}^+} \left[ c_{m,2p-1} \cos(px) + c_{m,2p} \sin(px) \right].
\end{equation}
 With the selection of suitable initial values and the application of relevant iterative optimization algorithms, a smooth Fourier representation $f_{m,\infty}(x)\in\mathcal{C}^{\infty}(0,2\Lambda)$ can be derived. We formally designate this overall approach as the \emph{Recursive Parametric Smoothing Strategy(RPSS)}. This name reflects its core mechanism: recursively refining the parameters (the Fourier coefficients) to find a smooth approximation of the true surface, thereby enabling an efficient inversion.


Having established the nonlinear equation \cref{eq:02} that governs the surface profile, we now introduce the numerical method for solving it. Our approach pairs the Recursive Parametric Smoothing Strategy (RPSS) with an iterative optimization algorithm to determine the unknown Fourier coefficients of the surface. For the iterative optimization, we adapt the Monte Carlo Continuation Method, an approach that has proven effective for similar inverse scattering problems. This method is well-suited for the nonlinearity of our problem, as the continuation strategy allows the algorithm to converge from a simple initial guess. Since the surface is uniquely determined by its Fourier coefficients $ \mathrm{c}_{m}=(c_{m,0},c_{m,1},\cdots,c_{m,n},\cdots)^{\mathsf{T}}$, the task of reconstructing the function is transformed into a finite-dimensional optimization problem for these coefficients. Therefore, the following sections will focus on the iterative optimization of this coefficient vector.


\subsection{The Monte Carlo Continuation method}{\label{MCCC}}
In practice, it is infeasible to take into account all the terms in the coefficient $\mathrm{c}_{m}$. Therefore we will pick a large truncation term number index $k_{max}$ and then perform the iterations in a finite-dimensional space. The reconstructed surface for the $m$-th realization, denoted by  $f_{m,k_{max}}(x)$ , is thus given by:
\begin{equation}{\label{fule}}
f_m(x):=f_{m,k_{max}}(x)=c_{m,0}+ \sum_{p\leqslant k_{max}} \left[ c_{m,2p-1} \cos(px) + c_{m,2p} \sin(px) \right],\end{equation} where $c_{m,2p-1}=c_{m,2p}=0\ \text{for all}\ p>k_{max}$ in \cref{fule}. Without loss of generality, define the spectral bandwidth index $k\ (k\leqslant k_{max})$ and coefficients vector 
\begin{equation}{\label{chen}}
    \mathrm{c}_{m,k}:=(c_{m,0},c_{m,1},\cdots,c_{m,2k},0,\cdots,0)^{\mathsf{T}}\in\mathbb{R}^{2k_{max}+1}
\end{equation}for the truncated coefficient vector. Our continuation method reconstructs these coefficients by progressively increasing the spectral bandwidth. Let $k$ (where $k\leqslant k_{max}$) be the truncation order at an intermediate step of the algorithm. For notational convenience, we can represent the corresponding coefficient vector \cref{chen} as an element in a infinite-dimensional space by padding it with zeros, i.e., $(c_{m,0},c_{m,1},\cdots,c_{m,2k_{max}},$ $0,0,\cdots)$. To simplify the notation throughout this paper, the vector $\mathrm{c}_{m,k}$ will be used to represent both the finite-dimensional vector in $\mathbb{R}^{2k_{max}+1}$ and its corresponding infinite-dimensional counterpart obtained by padding with zeros. This convention is unambiguous because our iterative algorithm, when operating at a bandwidth $k$, does not affect the coefficients of orders higher than $k$. 


The Monte Carlo Continuation method consists of two procedures. The initialization phase performs Monte Carlo sampling in $(\Omega,F,\mu)$ followed by a preliminary iteration, yielding a sparse coefficient vector whose non-zero entries are concentrated in the first $2k_1+1$ components $(k_1<k_{max})$,  with subsequent iterations progressively activating higher-indexed terms. Denote the total number of samples by $M$ and the prescribed wavenumber by $\kappa$. Let $k_{max}$ be the largest integer no larger than $\kappa$. For each sample $\omega_m=1,2,\cdots M$, the corresponding realization of the random surface is $f_m(x)=c_{m,0}+ \sum_{p=1}^{k_{max}} \left[ c_{m,2p-1} \cos(px) + c_{m,2p} \sin(px) \right]$, 
where the coefficients: 
\[\mathrm{c}_{m,k_{max}}=(c_{m,0},c_{m,1},\cdots,c_{m,2k_{max}})^{\mathsf{T}}\] are the parameters of the surface that is  to be updated recursively. Set the initial value $c_{m,0}=y_0$ and $ c_{m,p}=0, \ p=1,2,\cdots,2k_{max}$. Select $\kappa_1<\kappa$ as the initial wavenumber  and let $k=k_1$  be the largest integer no more than $\kappa_1$, then in this case  $\mathrm{c}$ will realize as :
\[\mathrm{c} := \mathrm{c}_{m,k_1} = \Bigl(  c_{m,0}, \dots, c_{m,2k_1},   \underbrace{0, \dots, 0}_{\mathclap{\substack{  2(k_{\max}-k_1)\text{ zero elements}}}} \Bigr)^{\!\mathsf{T}}.\]
Under these coefficients, the iterating surface will become:
\[f_{m,k_1}(x) = c_{m,0} + \sum_{p=1}^{k_1} \left[ c_{m,2p-1} \cos(px) + c_{m,2p} \sin(px) \right].\]
Considering the multi-angle incidence configuration with discrete directions $\theta_l\in(-\frac{\pi}{2},\frac{\pi}{2})$, $\ l=1,2,\cdots,L$, define
\begin{equation}{\label{553}}
      J_l(\mathrm{c}_{m,k_1}) = \left\| \sum_{n=-N}^{N} \psi_{n,l}(\omega_m) e^{i \alpha_n x + i \beta_n f(\omega_m;x)} + e^{i \alpha x - i \beta f(\omega_m,x)} \right\|^2_{\mathcal{L}^2(0,\Lambda)}. 
\end{equation}
Denote $\mathbf{J}_{k_1}:\mathbb{R}^{2k_1+1}\rightarrow\mathbb{R}^{L},\ \mathbf{J}_{k_1}(\mathrm{c}_{m,k_1}) = \left[J_1(\mathrm{c}_{m,k_1}),\cdots, J_L(\mathrm{c}_{m,k_1})\right]^{\mathsf{T}}
$ which collects all angles of incident wave information. Since \cref{eq:02}  holds for each incident angle $\theta_l$,  the composite problem can be formulated as:
\begin{equation}{\label{J}}
    \mathbf{J}_{k_1}(\mathrm{c}_{m,k_1}) = 0.
\end{equation}
 Complementing the previous iterative schemes like the Gauss-Newton approach \cite{bruckner2003two} and Levenberg-Marquardt method \cite{hettlich2002iterative}, the Landweber iteration provides an alternative first-order optimization scheme that achieves comparable accuracy with significantly reduced computational complexity. Let $\eta_k$ be a relaxation parameter dependent on the wavenumber, $\mathrm{c}_{m,k_1}^{(0)}=\mathrm{c}_{m,k_1}$ is the initial vector and the Jacobian matrix of \cref{J} is:
\[\mathbf{DJ}_{k_1} = \left( \frac{\partial J_l}{\partial c_{m,p}} \right)_{l=1,2,\ldots,L, p=0,1,\ldots,2k_1}.\]
The iterative expression for $c_{m,k}$ can then be written as follows:
\begin{equation}{\label{landwe}}
    \mathrm{c}_{m,k_1}^{(t+1)} = \mathrm{c}_{m,k_1}^{(t)} - \eta_\kappa \mathbf{D} \mathbf{J}_{k_1}^\mathrm{T}(\mathrm{c}_{m,k_1}^{(t)}) \mathbf{J}_{k_1}(\mathrm{c}_{m,k_1}^{(t)}), \quad t = 0, 1, 2, \ldots
\end{equation}
During the initialization phase, this recursive procedure will yield a coarse approximation after a few iterations, where the first $k_1$ non-zero coefficients capture the macroscopic  profile features while higher-order terms remain inactive. 
Up to this point, the first stage concludes.

In the subsequent completion phase, we will implement a continuation strategy to progressively capture high-frequency details by adaptively increasing the spectral bandwidth index $k$ from its initial value $k_1$ to the target resolution $k_{max}$. Select $\kappa_2>\kappa_1$ and pick $k_2>k_1$ the largest integer no larger than $\kappa_2$, then set $k=k_2$ in this round and follow the same strategy in the first stage. As for initialization, we let $\mathrm{c}_{m,k_2}=(c_{m,0},c_{m,1},\ldots,c_{m,2k_2})^{\mathsf{T}}$  be:
\[c_{m,p} := 
\begin{cases} 
c_{m,p}, & \text{for } 0 \leq p \leq 2k_1, \\
0, & \text{for } 2k_1 < p \leq 2k_2,
\end{cases}\]
and $\mathrm{c}^{(0)}_{m,k_2}=\mathrm{c}_{m,k_2}$. The continuation stage terminates when $\kappa_1,\kappa_2,\ldots$ reaches the threshold wavenumber $\kappa$ or meets the stopping criterion. Upon completion of the MCC iteration, each sample $\omega_m$ is associated with an optimal set of coefficient parameters $\mathrm{c}_{m,k_{max}}=(c_{m,0},c_{m,1},\cdots,c_{m,2k_{max}})^{\mathsf{T}}$ derived through this recursive process.

\subsection{Modified Monte Carlo continuation method}{\label{MC}}
As demonstrated in \cite{gu2024inverse}, distinct realizations of the same stochastic structure exhibit consistent macroscopic properties,  leveraging this invariance to initialize the subsequent iterative samples will enable significant computational savings. In other words, MCC requires only a limited number of samples to determine the surface's macroscopic properties. Therefore, the rough outlines obtained from these initial samples can serve as substitutes for the first several iterations of other samples in subsequent MCC calculations. More specifically, the first $M_r(M_r<M)$ samples will undergo the complete iteration process in \Cref{MCCC}, the rest samples will start their iteration from  initialization $f(x)=\frac{1}{M_r}\sum_{m=1}^{M_r}f_m(x)$ directly. Following the effectiveness of this method in  \cite{wang2025numerical}, we  incorporate this technique to optimize computational efficiency.

\subsection{The inversion of statistical functions $g(x)$ and $h^2(x)$\ (MCCh)}{\label{acc}}
Since \[\mathbb{E}[f(\omega,x)] = g(x) + \underbrace{\mathbb{E}\left[\sum_{j=0}^{N_0} h(x_j) \xi_j(\omega) \sqrt{\Delta x} \varphi_j(x)\right]}_{= 0}= g(x),
\]we first reconstruct  the mean profile $g(x)$ from our samples. Analogous to the use of empirical risk minimization in machine learning as a surrogate for expected risk minimization, we recover the dominant profile of the original stochastic surface $g(x)$ by by computing the sample average:
\[g(x)\approx\bar{f}(x)=\frac{1}{M}\sum_{m=1}^{M}f_m(x).\]
Equivalently, this can be expressed in terms of the coefficients $\mathrm{c}_{m,k_{max}}$ for each realization: 
\begin{equation}{\label{334}}
    \mathbf{c}=\frac{1}{M}\sum_{m=1}^{M}\mathrm{c}_{m,k_{max}}=\frac{1}{M}\Big(\sum_{m=1}^{M}c_{m,0},\sum_{m=1}^{M}c_{m,1},\cdots,\sum_{m=1}^{M}c_{m,2k_{max}}\Big).
\end{equation}
Thus, the  the mean profile $g(x)$ of the random periodic structure can be reconstructed via a truncated Fourier series:
\[g(x)\approx\bar{f}(x) = \bar{c}_0 + \sum_{p=1}^{k_{max}} [\bar{c}_{2p-1} \cos(px) + \bar{c}_{2p} \sin(px)].\]
The inversion of  the variance intensity  $h^2(x)$ requires computation of the covariance matrix. Recalling that for a stochastic process, the covariance function $c(\cdot,\cdot)$ is defined as:
\[c(s,t) = \mathbb{E}\left[f(s) - \mathbb{E}f(s)\right]\left[f(t) - \mathbb{E}f(t)\right], \quad s,t \in (0,2\pi).\]
In our framework, we approximate this covariance function using the empirical estimator:
\[c(s,t) = \frac{1}{M} \sum_{m=1}^M [f_m(s) - \bar{f}(s)] [f_m(t) - \bar{f}(t)], \quad s, t \in (0, 2\pi),\]
where $\bar{f}(x)$ represents the empirical mean function defined previously.
To ensure unbiased estimation, the coefficient in the covariance formula should theoretically be $\frac{1}{M-1}$. However, since our sample size $M$ is sufficiently large in practice ($M \gg 1$), we adopt $\frac{1}{M}$ as the normalization factor for notational convenience. For spatial discretization, we maintain consistency with the forward problem's scheme. Let $\Delta x = x_i - x_{i-1}$ denote the uniform grid spacing over the interval $[0, 2\pi]$, where ${x_i}$ represents the discrete nodal points.  We can deduce that:
\begin{align}{\label{hin}}
c_{ij}=cov(x_i,x_j)&=\mathbb{E}\Bigg( \Big(f_m(x_i)-\mathbb{E}f(x_i)\Big)\Big(f_m(x_j)-\mathbb{E}f(x_j)\Big)\Bigg)\\&=\mathbb{E}\Bigg(\Big(\sum_{k=0}^{N} h(x_k) \xi_k(\omega) \sqrt{\Delta x} \varphi_k(x_i)\Big)\Big(\sum_{k=0}^{N} h(x_k) \xi_k(\omega) \sqrt{\Delta x} \varphi_k(x_j)\Big)\Bigg)\notag\\&=\mathbb{E}\Bigg(\Big(h(x_i)\xi_i\sqrt{\Delta x}\Big)\Big(h(x_j)\xi_j\sqrt{\Delta x}\Big)\Bigg)\notag
\\&=h(x_i)h(x_j)\Delta x\mathbb{E}(\xi_i,\xi_j)\quad\quad\text{where}\  \xi_i,\xi_j\thicksim\mathcal{N}(0,1)\notag.
\end{align}
The equation is non-trivial only when $i= j$, thus we can get:
\begin{equation}{\label{aa}}
    c_{ii}=h^2(x_i)\Delta x.
\end{equation}
Then we can reconstruct  the variance intensity $h^2(x)$ valued at $x_i$ as:
\[h^2(x_i)=\frac{c_{ii}}{\Delta x}.\]
Our analysis leads to the algorithm in \cref{MCCh}.

\begin{algorithm}[!htbp] 
\caption{MCCh method}
\label{MCCh}
\begin{algorithmic}[1]  
\STATE Input grid spacing hyperparameters  $\Delta x$, Maximum wavenumber $k_{max}$
\STATE{Define $\mathbf{c}$ as in \cref{334} and set $\mathbf{c}=0$}
\STATE{Determine the sample size $M_r$ for the initial phase of iteration}
\STATE{choose $y_0$, define the maximum allowable iteration times $T$}

\textbf{Step 1: Initialization Stage} 
\FOR{$m = 1$ \TO $M_r$}{
  \STATE Generate $\omega_m$ and define  $\mathrm{c}_{m,k_{max}}=(c_{m,0},c_{m,1},\cdots,c_{m,2k_{max}})^{\mathsf{T}}$
  \STATE Set $c_{m,0}=y_0, c_{m,p}=0,p=1,2,\cdots,2k_{max}$
  \FOR{$k=k_1,k_2,\cdots,k_{r}$}{
  \STATE Define $\mathrm{c}_{m,k}^{(0)}=(c_{m,0},c_{m,1},\cdots,c_{m,2k})^{\mathsf{T}}$
 \FOR{$t=0,1,2,\cdots,T$ }{ 
   \STATE Calculate $\alpha,\beta$ for each incident angle $\theta_l\in(-\frac{\pi}{2},\frac{\pi}{2})$, $l=1,2,\cdots,L$ }
   \STATE Define $\mathbf{J}_{k}(\mathrm{c}_{m,k}) = \left[J_1(\mathrm{c}_{m,k}),\cdots, J_L(\mathrm{c}_{m,k})\right]^T$ where $J_i$ is defined as in \cref{553}
   \STATE Calculate Jacobian matrix $\mathbf{DJ}_{k} = \left( \frac{\partial J_l}{\partial c_{m,p}} \right)_{l=1,2,\ldots,L, p=0,1,\ldots,2k}$
   \STATE $\mathrm{c}_{m,k}^{(t+1)} = \mathrm{c}_{m,k}^{(t)} - \eta_k \mathbf{D} \mathbf{J}_{k}^\mathrm{T}(\mathrm{c}_{m,k}^{(t)}) \mathbf{J}_{k}(\mathrm{c}_{m,k}^{(t)})$
  \ENDFOR
  \STATE $c_{m,p} := c_{m,p}^{(T+1)}, \ p = 0, 1, \ldots, 2k$
  }
  \ENDFOR
  \STATE Let $f_m(x) = c_{m,0} + \sum_{p=1}^{k_{\text{max}}} \left( c_{m,2p-1} \cos px + c_{m,2p} \sin px \right)$
  }
\ENDFOR  
\STATE Let $\hat{\mathrm{c}}=\frac{1}{M_r}\sum_{m=1}^{M_r}\mathrm{c}_{m,k_{max}}$

~\\
\textbf{Step 2: Completion Stage} 
\FOR{For $m=M_r+1,M_r+2,\cdots,M$}{
\STATE Generate $\omega_m$ and define  $\mathrm{c}_{m,k_{max}}=(c_{m,0},c_{m,1},\cdots,c_{m,2k_{max}})^{\mathsf{T}}$
\STATE Set $\mathrm{c}_{m,p} = \hat{\mathrm{c}}$
  \FOR{$k=k_r,k_{r+1},\cdots,k_{max}$}{
  \STATE Repeat the same recursive process in Step 1 with larger $k$
  }
  \ENDFOR
}
\ENDFOR

 ~\\
\textbf{Step 3: Inversion Stage}
\STATE Calculate $\mathbf{c}=\frac{1}{M} \sum_{m=1}^{M} c_{m, k_{\text{max}}}$
\STATE Let  $\bar{f}(x) = \bar{c}_0 + \sum_{p=1}^{k_{\text{max}}} \left[ \bar{c}_{2p-1} \cos(px) + \bar{c}_{2p} \sin(px) \right]$
\STATE Calculate the covariance matrix $\mathbf{C}=(c_{ij})$, 
\STATE where $c_{ij}=cov(x_i,x_j)=\frac{1}{M} \sum_{m=1}^M\Big( \big(f_m(x_i)-\bar{f}(x_j)\big)\big(f_m(x_i)-\bar{f}(x_j)\big)\Big)$
\STATE Extract the diagonal elements of $\mathbf{C}$ as $\text{diag}(\mathbf{C}):=(c_{11},c_{22},\cdots,c_{nn})$ and denote $\mathbf{x}:=(x_0,x_1,\cdots,x_n)$
\STATE Calculate $h^2(\mathbf{x})=\frac{\text{diag}(\mathbf{C})}{\Delta x}$ or $|h(\mathbf{x})|=\sqrt{\frac{\text{diag}(\mathbf{C})}{\Delta x}}=\sqrt{\frac{1}{\Delta x}}(c_{11},c_{22},\cdots,c_{nn})$ \\
(If prior geometric information about $h(x)$ is available, further reconstruction of $h(x)$ may be possible using statistical techniques.)
\RETURN $\bar{f}$ and $h^2(\mathbf{x})$ or $|h(\mathbf{x})|$  
\end{algorithmic}
\end{algorithm}

\begin{remark}[$h^2(x)$ VS $h^2(\mathbf{x})$]
   Following the modeling in \Cref{sec:For}, the white noise increments are sampled exclusively at the grid points ${x_i}$. Consequently, the estimated values ${h^2(x_i)}$ derived from \cref{aa} actually represent discrete approximations of the local variance intensity $h^2(x)$, capturing the sample variance intensity of $h^2(x)$ only at these nodal points. Therefore, strictly speaking, the output of our algorithm is the vector $h^2(\mathbf{x})$, where $\mathbf{x}=(x_0,x_1,\ldots,x_{N_0})$ represents the grid points, and $h^2(\mathbf{x})$ is the vector of estimated variance intensities at these points. Provided that $|\Delta x|$ is sufficiently small and satisfies the consistency requirements for both forward and inverse computations, this discrete approximation will converge to the continuous function within an acceptable error tolerance. Although refining the grid could theoretically enhance the inversion accuracy for $h^2(x)$, such refinement would inevitably introduce stronger random oscillations and significantly increase computational costs. Therefore, selecting an optimal $\Delta x$ that balances numerical stability with reconstruction precision becomes crucial in practice. Given this convergence for a sufficiently small $\Delta x$, we will adopt the continuous notation $h(x)$ (or $|h(x)|$) to refer to our discrete vector reconstruction $h(\mathbf{x})$ (or $|h(\mathbf{x})|$) in the subsequent discussions, with the understanding that it represents a pointwise approximation of the true function.
\end{remark}

\begin{remark}[Possibility of reconstructing $h(x)$]
    As we can see, our method ultimately yields $g(x)$ and $h^2(x)$. The inversion of  $h(x)$ is hindered by our lack of knowledge regarding the true sign of $h(x)$. Moving forward, if we manage to obtain some prior information about $h(x)$, such as its sign and shape information, and combine this with widely used statistical approaches for data analysis\begin{tikzpicture}[baseline=-0.65ex]
    \draw (0,0) -- (0.5,0); 
\end{tikzpicture}extracting the commonalities among majority of samples meanwhile discard a few outliers\begin{tikzpicture}[baseline=-0.65ex]
    \draw (0,0) -- (0.5,0); 
\end{tikzpicture}we will achieve an accurate inversion of $h(x)$. Additional details are to be presented in the discussion parts \Cref{hxinv} later.
\end{remark}

\section{Experimental results}{\label{sec:experiments}}

In this part we will present the numerical  experiments to verify the effectiveness of  our inversion algorithm. At start, we  need to first acquire the near-field measured scattering data $u(\omega_m;x,y_0)$ from the  forward problem with each surface $f_{lip}^m(x)$. Here we will adopt the adaptive finite element method with perfectly matched absorbing layers \cite{chen2003adaptive} to simulate $u(\omega_m;x,y_0)$. To ensure that all random perturbations in the profile are accurately captured during geometric meshing process, we incorporate $\Delta x$ as a key parameter in the forward problem solution.  As the grid is refined with decreasing $\Delta x$, the geometric mesh scale diminishes accordingly, leading to increased computational burden. The introduction of \cref{110} effectively constrains our computational overhead within manageable limits. Besides, in order to emulate the actual situation, data with noise will be utilized as the true measured scattered data. Formally, the scattering data will have the form:
\begin{equation}{\label{noi}}
    u(\omega; x, y_0) := u(\omega; x, y_0)(1 + \tau \epsilon).
\end{equation}
The term $\epsilon$ indicates a set of random numbers uniformly distributed over the interval $[-1,1]$, and $\tau$ is the factor that quantifies the noise level in the data measurements, in the following examples we will take $\tau=0.1\%$. Typically, the Monte Carlo approach necessitates the collection of a substantial number of samples. In the numerical examples, we find that for our given cases, a sample size of around $M=10^3$ usually suffices to attain an acceptable level of precision. Besides, for the starting coefficient vector $\hat{\mathrm{c}}$ in the Completion Stage in  \cref{MCCh}, a sample size of around $M_r=10^2$ from the first phase is sufficient to depict the general macroscopic profile of the surface. $\eta_k$ in \cref{landwe} will be set as $\eta_\kappa=10^{-3}/\kappa$ at each stage and truncation number $N$ in \cref{eq:02} will take $N=8$. To avoid overly dense oscillations which would violate the constraints of \cref{110}, we set the maximum number of partition points on to $N_0=160$, and this results in a step size of  $\Delta x\approx 0.004$. Additionally, to maintain periodicity condition $f(0)=f(2\pi)$, we enforce the noise values at the boundary points $x_0$ and $x_n$ to be identical through constrained randomization: $\xi(x_0) \equiv \xi(x_n)$. 

Within the framework of our algorithm, the precise inversion of $g(x)$ is particularly crucial, as it significantly determines the accuracy of the subsequent inversion of $h^2(x)$. 
In the following numerical examples, we will observe that, under the premise that there is a sufficient number of samples to support the reconstruction and the shape of $g(x)$ can be well reconstructed, our method can accurately invert $|h(x)|$ when $h(x)$ belongs to a simple linear combination of bounded elementary functions.

We first begin with a simple example:

\begin{example}{\label{ex1}}
     Reconstruct the random surface where
\[g_1(x)=0,\quad h_1(x)=\cos(x),\] inside the interval $[0,2\pi]$. 
\end{example}  

\begin{figure}[htbp]
    \centering
    \begin{subfigure}{0.327\textwidth}
        \includegraphics[width=\linewidth]{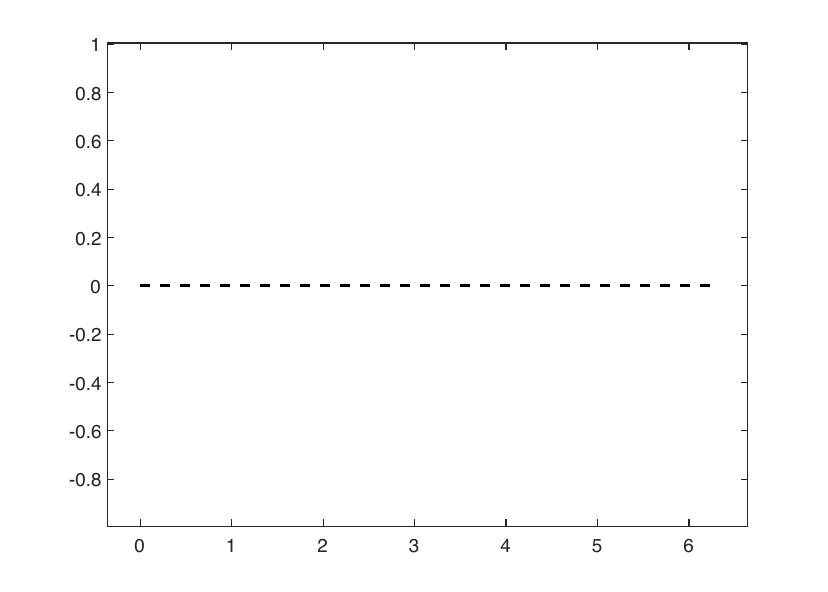}
        \caption{\it{\small Noise-free surface $g_1(x)$ }}
        \label{null}
    \end{subfigure}
     \begin{subfigure}{0.325\textwidth}
        \includegraphics[width=\linewidth]{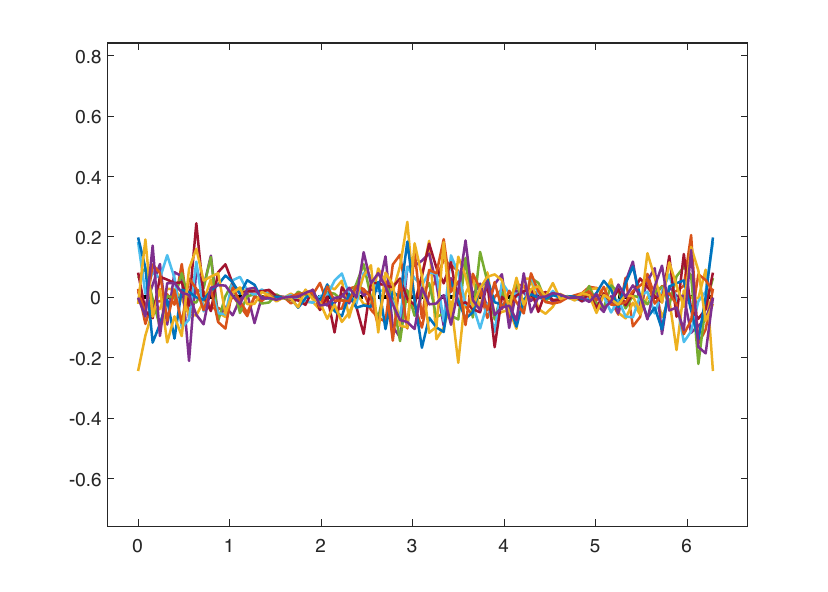}
        \caption{\it{\small Different realizations }}
        \label{e0}
    \end{subfigure}
    \begin{subfigure}{0.325\textwidth}
        \includegraphics[width=\linewidth]{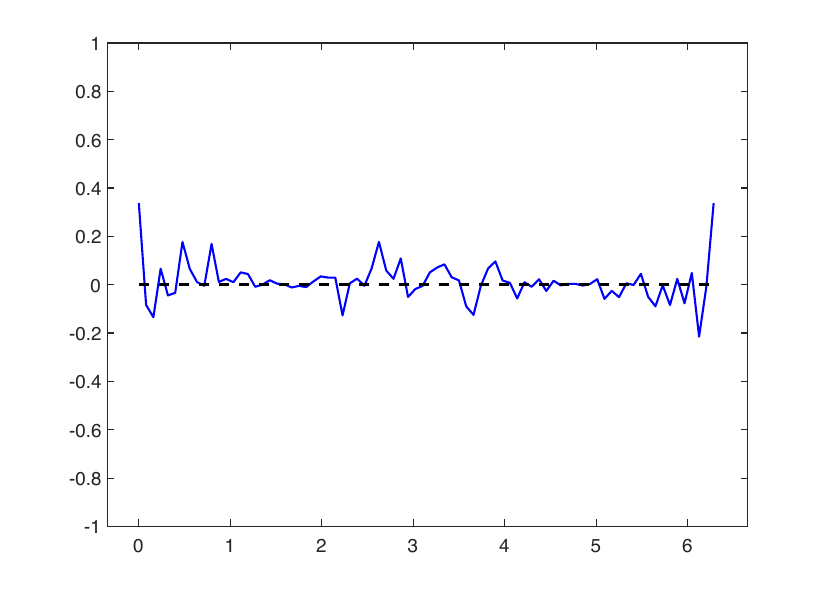}
        \caption{\it{\small The $10^{th}$ surface  $f^{10}_{lip}(x)$}}
        \label{10e0}
    \end{subfigure}
     \caption{\it{\small Numerical simulation of \cref{ex1} at spacial resolution $N_0=80$. (a): The image of noise-free function $g_1(x)$. (b): The realizations of first 10 samples $\{f^m_{lip}(x)\}_{m=1}^{10}$. (c): The tenth concrete realization $f^{10}_{lip}(x)$ of random surface.}}
    \label{fig0}
\end{figure}

\begin{figure}[htbp]
    \centering
    \begin{subfigure}{0.325\textwidth}
        \includegraphics[width=\linewidth]{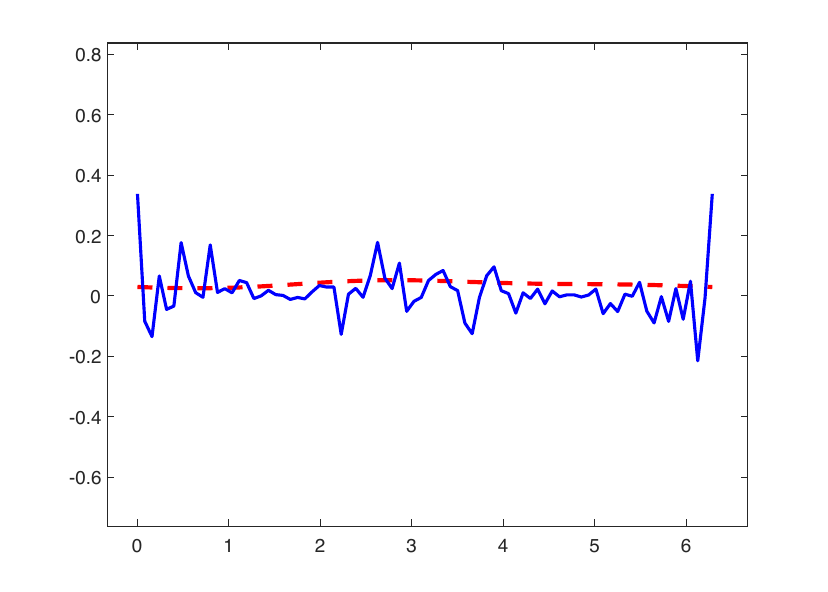}
         \caption{\it{\small Reconstructed   $f_{10}(x)$}}
        \label{heyu}
    \end{subfigure}
     \begin{subfigure}{0.325\textwidth}
        \includegraphics[width=\linewidth]{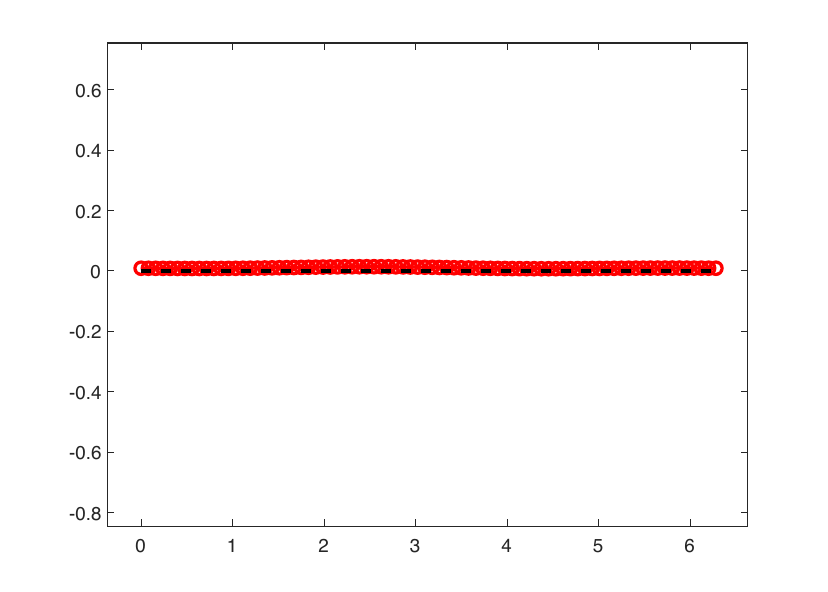}
        \caption{\it{\small Reconstructed $\bar{f}(x)$ } }
        \label{mayu}
    \end{subfigure}
    \begin{subfigure}{0.325\textwidth}
        \includegraphics[width=\linewidth]{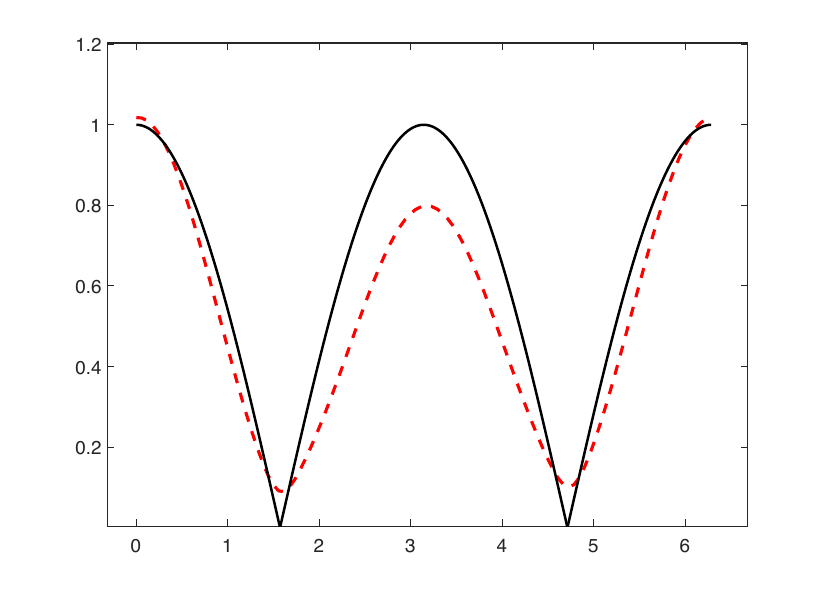}
        \caption{\it{\small Reconstructed $|h_1(x)|$}}
        \label{acccc}
    \end{subfigure}
    \caption{\it{\small Inversion of $g_1(x)$ and $h_1(x)$. (a): Solid  line: 10th surface $f^{10}_{lip}(x)$. Dashed  line: the inversion result $f_{10}(x)$ obtained at the end of the tenth algorithm execution. (b): Dashed line: the original image of $g_1(x)$. Circular marked line: the reconstructed $\bar{f}(x)$ . (c):The inversion of $|h_1(x)|$ in \cref{ex1}. Solid  line: the original image of $|\cos(x)|$. Dashed  line: the  image of reconstructing $|h_1(x)|$.}}  
 
\end{figure}

We initialize $N_0=80$ which means $\Delta x\approx 0.079$.  The left panel \cref{null} in \cref{fig0} shows the image $g_1(x)$ without any random noise, whereas the right-hand side image \cref{10e0} displays the  specific realization of the perturbed random surface at the 10th iteration. First $10$ realizations of random surfaces are presented in the middle image \cref{e0}. Since $g_1(x)$ is a constant function in \Cref{ex1}, we observe that the inversion of $g_1(x)$ performs best when $\kappa=2$. Set $k_{max}=2$. The maximum number of Landweber iterations is set to be 50. Here to exemplify the inversion performance, we present the results from the 10th algorithmic iteration applied to the 10th sample in \cref{heyu}.

\cref{mayu} shows the image of reconstructed $\bar{f}(x)$. 
By comparing the reconstructed coefficients with their theoretical values in \cref{table0}, the inversion results of $g_1(x)$ exhibit satisfactory accuracy within an acceptable error margin.  Since $g_1(x)=0$ in \Cref{ex1}, it is treated as a Fourier series characterized by a null set of coefficients, as listed in the `Target' row in \Cref{table0}. 

\begin{table}[htbp]
  \centering
  \begin{tabular}{lccccccc}
    \toprule
    & $c_{m,0}$ & $c_{m,1}$ & $c_{m,2}$ & $c_{m,3}$ & $c_{m,4}$  \\
    \midrule
    Target\ \big($g_1(x)$\big)    & 0   & 0 & 0 & 0 & 0        \\
    Initial   & 0.5617   & 0     & 0     & 0     & 0         \\
    Computed\ \big($\bar{f}(x)$\big)  & $0.0060$   & $-0.0016$ & $0.0016$ & 0.00064 & $-0.0017$  \\
    \bottomrule
  \end{tabular}
   \caption{\it{\small Accuracy of reconstructed Fourier coefficients for $g_1(x)$ versus reconstructing $\bar{f}(x)=c_{m,0} + \sum_{p=1}^{2} [ c_{m,2p-1} \cos(px) + c_{m,2p} \sin(px) ]$ of \cref{ex1}. Target $g_1(x)$: Theoretical Fourier coefficients of the ground truth function. Initial: Initial guess for coefficients. Computed $\bar{f}(x)$: Reconstructed mean estimate   with $k_{max}=2$.}}  
  \label{table0}
\end{table}

\begin{figure}[htbp]
    \centering
    \begin{subfigure}{0.325\textwidth}
        \includegraphics[width=\linewidth]{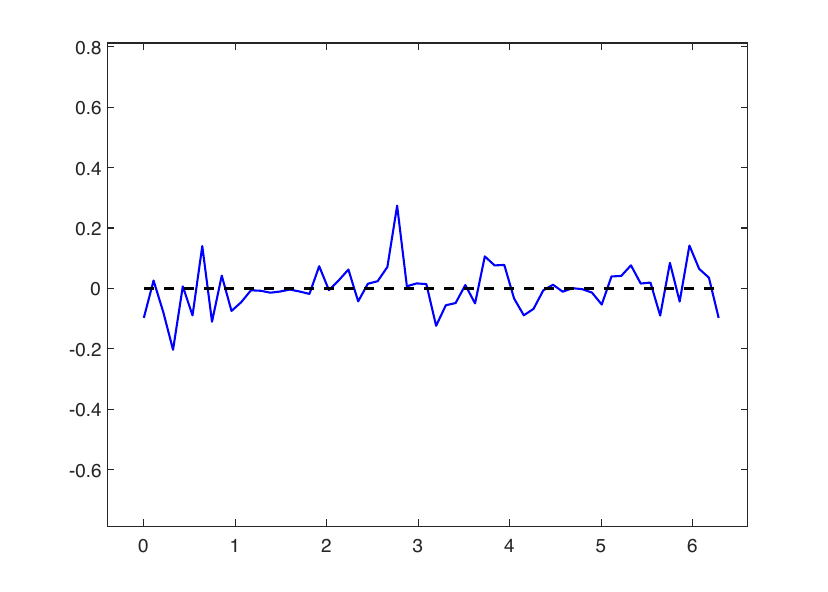}
        \caption{\it{\small   $f^{10}_{lip}(x)$ with $N_0=60$}}
        \label{180}
    \end{subfigure}
    \begin{subfigure}{0.325\textwidth}
        \includegraphics[width=\linewidth]{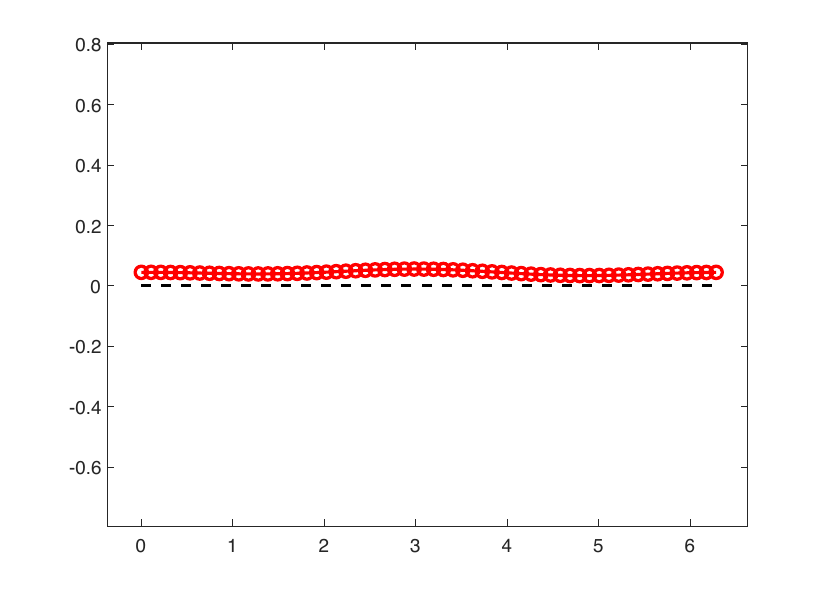}
        \caption{\it{\small  $\bar{f}(x)$ with $N_0=60$}}
        \label{180ginv}
    \end{subfigure}
      \begin{subfigure}{0.325\textwidth}
        \includegraphics[width=\linewidth]{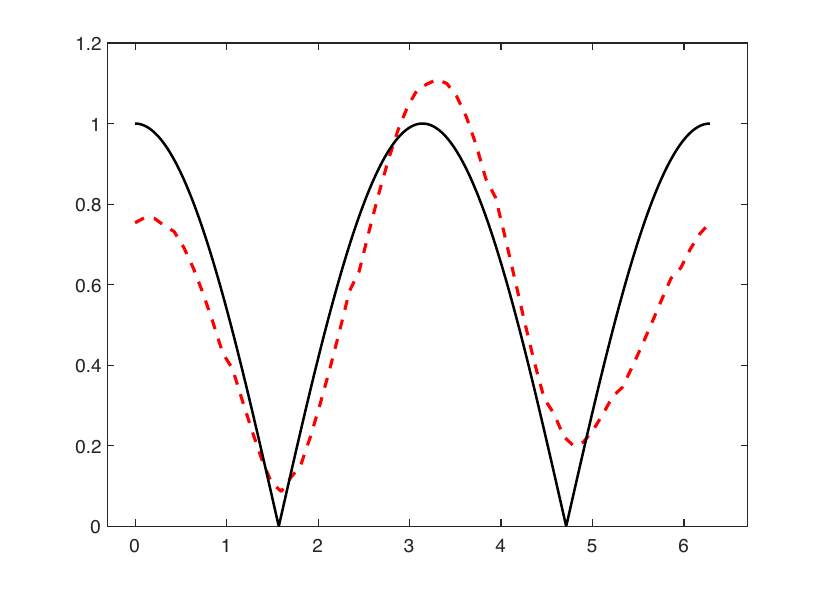}
        \caption{\it{\small   $|h_1(x)|$ with $N_0=60$}}
        \label{180inv}
    \end{subfigure}

     \begin{subfigure}{0.325\textwidth}
        \includegraphics[width=\linewidth]{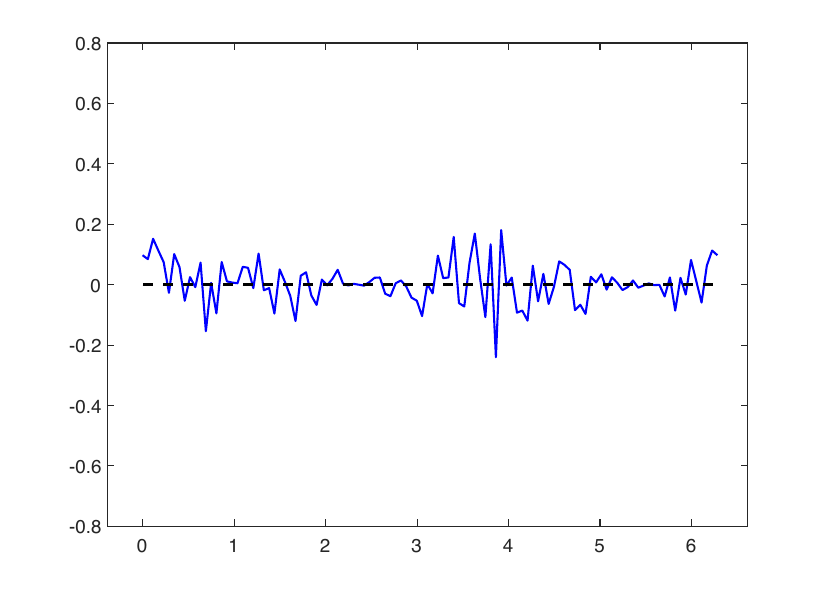}
        \caption{\it{\small $f^{10}_{lip}(x)$ with $N_0=110$ } }
        \label{1180}
    \end{subfigure}
\begin{subfigure}{0.325\textwidth}
        \includegraphics[width=\linewidth]{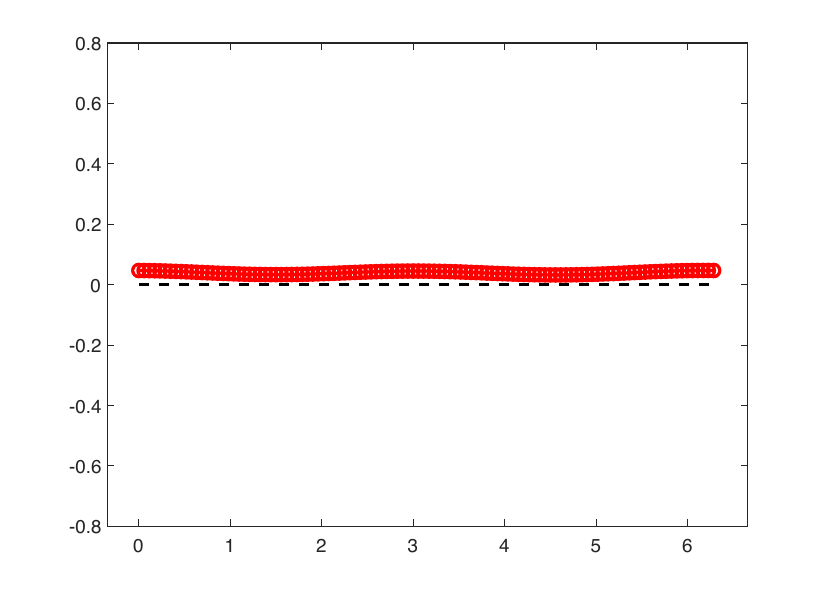}
        \caption{\it{\small  $\bar{f}(x)$ with $N_0=110$ } }
        \label{1180ginv}
    \end{subfigure}
    \begin{subfigure}{0.325\textwidth}
        \includegraphics[width=\linewidth]{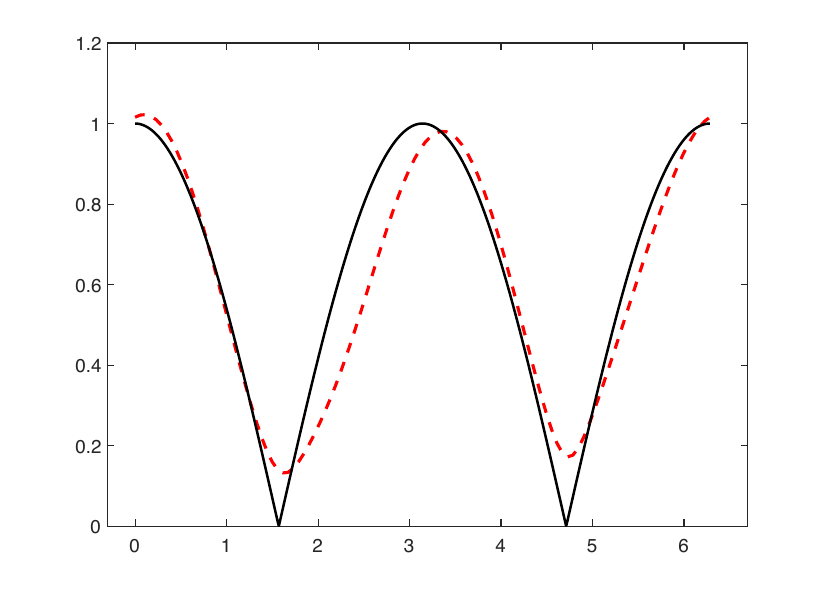}
        \caption{\it{\small  $|h_1(x)|$ with $N_0=110$ } }
        \label{1180inv}
    \end{subfigure}
    \caption{\it{\small The realization and the inversion at two different spatial resolutions for \Cref{ex1}. The top row corresponds to a coarse grid $N_0=60$ while the bottom row uses a finer grid $N_0=110$. First column:  the 10-th realization of random surface at two different spatial resolutions (Solid line) and  original $g_1(x)$ (Dashed line);  Second column: the reconstructed mean profile $\bar{f}(x)$ (Circular markers) compared with the ground truth $g_1(x)$(Dashed line).  Third column: the reconstructed variance intensity $|h_1(x)|$ (Dashed line)  compared with the true function $|\cos(x)|$ (Solid line).}}
    \label{ex0inv}
    \end{figure}

The inversion of $|h_1(x)|$ is shown in \cref{acccc}. Subsequently, we further examined the cases when $\Delta x=60$ and $\Delta x=110$ respectively. The corresponding numerical results are presented in \cref{ex0inv}. These pictures  demonstrate  that our algorithm is capable of precisely reconstructing $|h_1(x)|$ provided that the assumptions outlined in \Cref{110} and \Cref{911} regarding $\Delta x$ are satisfied. Both reconstructions are achieved within an acceptable error margin, thereby validating our algorithm's reliability under the specified conditions.


\begin{example}{\label{ex2}}
    Reconstruct the random surface where
\[g_2(x)=1.5+0.2\cos(x)+0.2\cos(2x),\quad h_2(x)=\sin(x),\] inside the interval $[0,2\pi]$.
\end{example} 

\begin{figure}[htbp]
    \centering
    \begin{subfigure}{0.327\textwidth}
        \includegraphics[width=\linewidth]{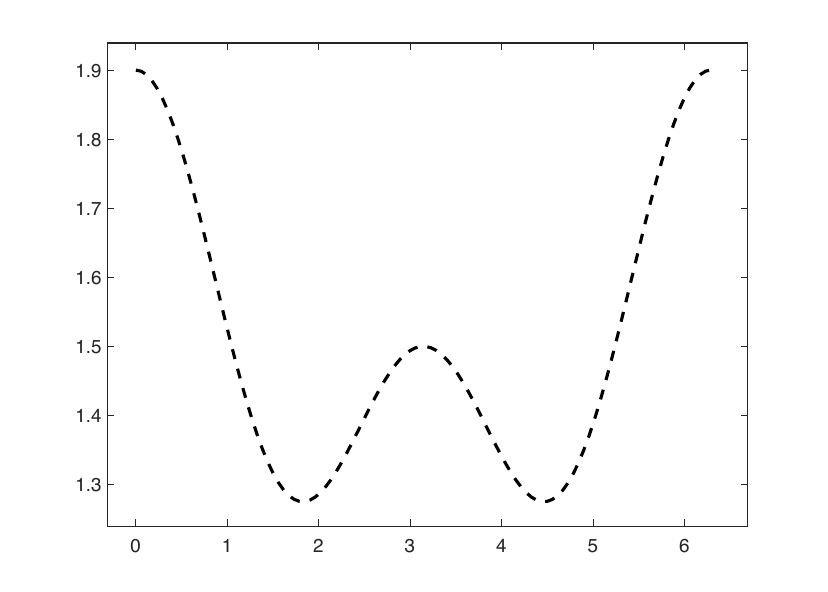}
        \caption{\it{\small Noise-free surface $g_2(x)$ }}
        \label{fig:image1}
    \end{subfigure}
     \begin{subfigure}{0.325\textwidth}
        \includegraphics[width=\linewidth]{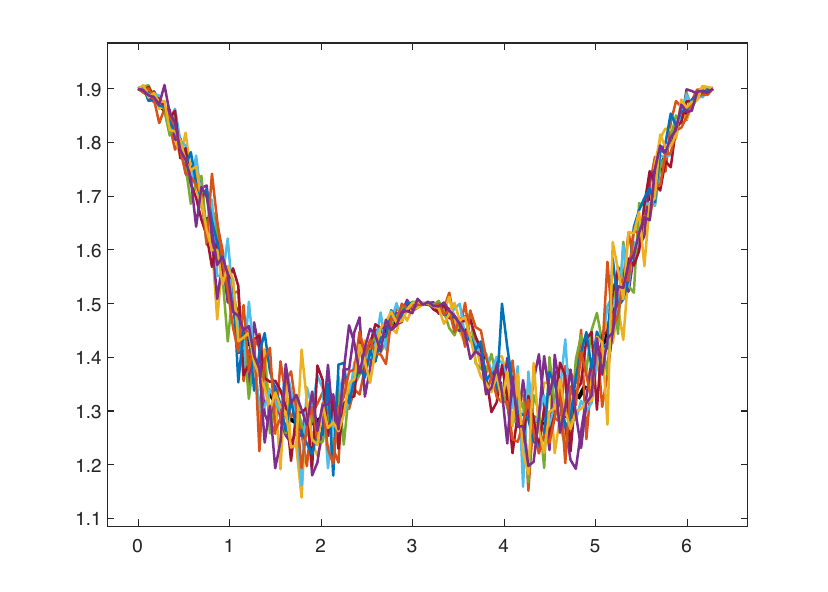}
        \caption{\it{\small First 10 realizations }}
        \label{fig:image2}
    \end{subfigure}
    \begin{subfigure}{0.325\textwidth}
        \includegraphics[width=\linewidth]{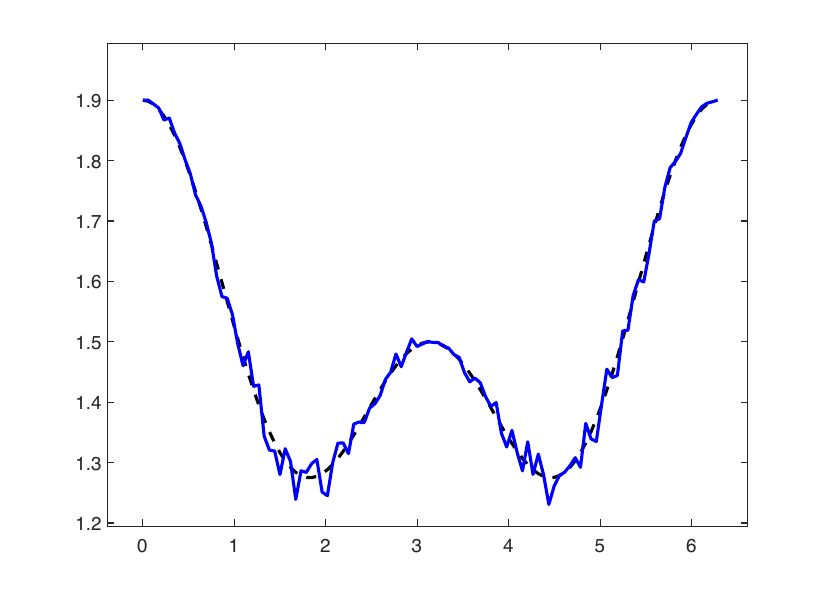}
        \caption{\it{\small The $10^{th}$ surface  $f^{10}_{lip}(x)$}}
        \label{fig:image3}
    \end{subfigure}
     \caption{\it{\small Numerical simulation of \cref{ex2} at spacial resolution $N_0=110$. (a): The image of noise-free function $g_2(x)$. (b): The realizations of first 10 samples $\{f^m_{lip}(x)\}_{m=1}^{10}$. (c): The tenth concrete realization $f^{10}_{lip}(x)$ of random surface.}}
    \label{figg}
\end{figure}

For this  example,  we first examine the intermediate discretization with 
$N_0=110$ nodes, corresponding to $\Delta x\approx 0.057$. The left panel \cref{fig:image1} in \cref{figg} illustrates the image $g_2(x)$ without any random noise,  the right-hand side image \cref{fig:image3} displays the  specific realization of the perturbed random surface at the 10th iteration. First $10$ realizations of random surfaces are presented in the middle image \cref{fig:image2}.

\begin{figure}[htbp]
    \centering
    \begin{subfigure}{0.325\textwidth}
        \includegraphics[width=\linewidth]{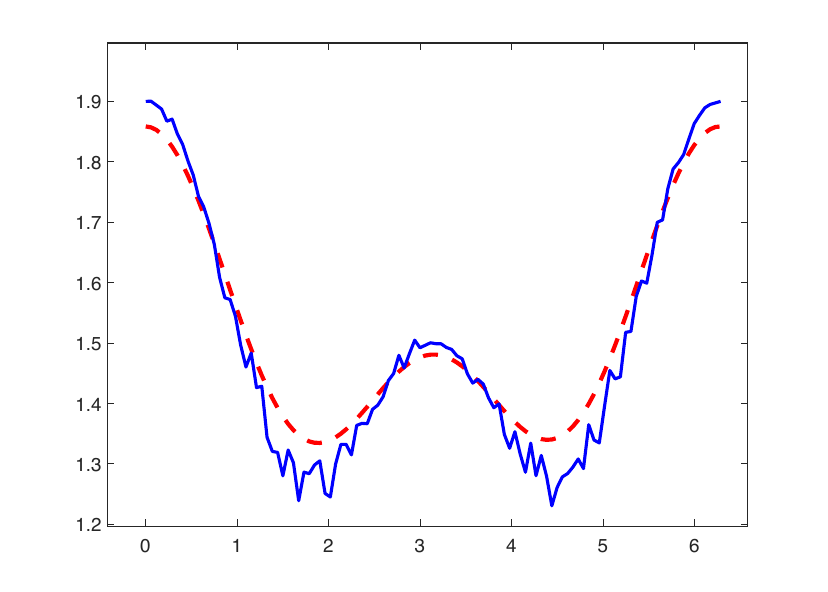}
        \caption{\it{\small Reconstruction $f_{10}(x)$}}
        \label{121212}
    \end{subfigure}
     \begin{subfigure}{0.325\textwidth}
        \includegraphics[width=\linewidth]{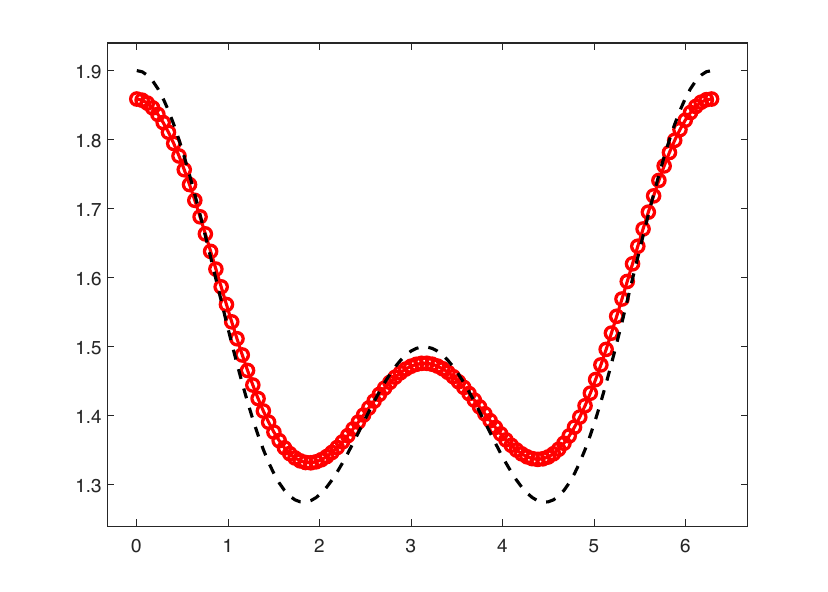}
        \caption{\it{\small Reconstructed $\bar{f}(x)$ } }
        \label{556}
    \end{subfigure}
     \begin{subfigure}{0.325\textwidth}
  \includegraphics[width=\linewidth]{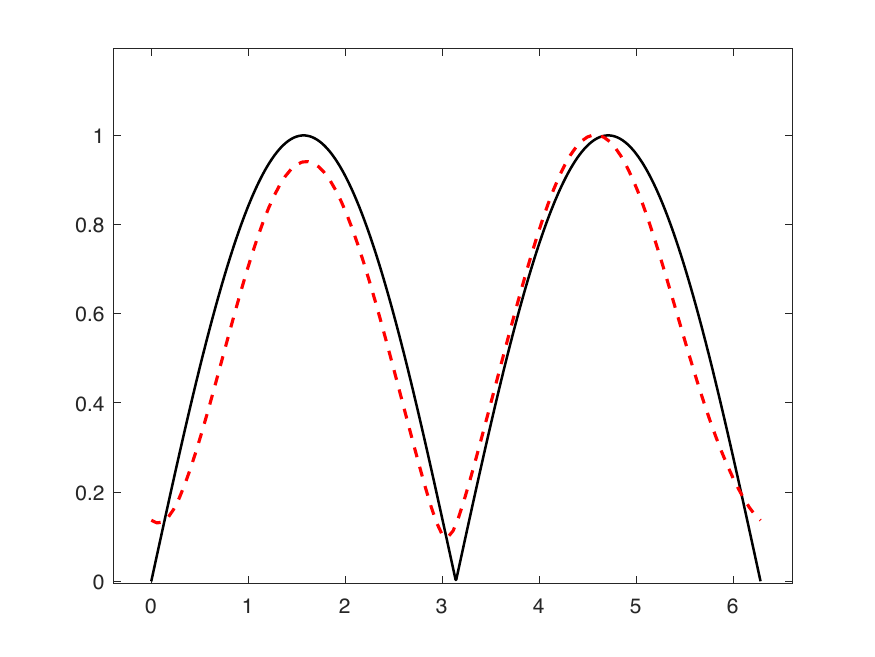}
  \caption{Reconstructed $|h_2(x)|$}
        \label{13}
    \end{subfigure}
    \caption{\it{\small Inversion of $g_2(x)$ and $h_2(x)$. (a): Solid  line: 10th surface $f^{10}_{lip}(x)$. Dashed  line: the inversion result $f_{10}(x)$ obtained at the end of the tenth algorithm execution. (b): Dashed line: the original image of $g_2(x)$. Circular marked line: the reconstructing $\bar{f}(x)$ . (c): The inversion of $|h_2(x)|$ in \cref{ex2}. Solid  line: the original image of $|\sin(x)|$. Dashed  line: the  image of reconstructing $|h_2(x)|$.}}  
    \end{figure}

Through parametric comparison, the data combination $\kappa=2$ and $k_{max}=2$  was identified as providing superior inversion accuracy compared to other tested configurations. The reconstructed profile  $f_{10}(x)$  visualized in \cref{121212} (compared with  $f^{10}_{lip}(x)$ ) is  obtained through consistent implementation of the Landweber scheme described in  \Cref{ex1}.  
From the acquired dataset, we obtained ultimate $\bar{f}(x)$ depicted in \cref{556}. As evidenced by the tabulated data \cref{tabel1}, the inverted $\bar{f}(x)$ exhibits close agreement with the ground truth $g_2(x)$.
\begin{table}[htbp]
  \centering
  \begin{tabular}{lccccccc}
    \toprule
    & $c_{m,0}$ & $c_{m,1}$ & $c_{m,2}$ & $c_{m,3}$ & $c_{m,4}$  \\
    \midrule
    Target\ \big($g_2(x)$\big)    & 1.5   & 0.2 & 0 & 0.2 & 0        \\
    Initial   & 2.1003   & 0     & 0     & 0     & 0         \\
    Computed\ \big($\bar{f}(x)$\big)  & $1.5162$   & $0.1913$ & $-0.0034$ & 0.1514 & $-0.0012$  \\
    \bottomrule
  \end{tabular}
   \caption{\it{\small Accuracy of reconstructed Fourier coefficients for $g_2(x)$ versus reconstructing $\bar{f}(x)=c_{m,0} + \sum_{p=1}^{2} [ c_{m,2p-1} \cos(px) + c_{m,2p} \sin(px) ]$ of \cref{ex2}. Target $g_2(x)$: Theoretical Fourier coefficients of the ground truth function. Initial: Initial guess for coefficients. Computed $\bar{f}(x)$: Reconstructed mean estimate  with $k_{max}=2$.}} 
  \label{tabel1}
\end{table}
Following  \cref{MCCh}, we obtained the inversion results of  $|h_2(x)|$ presented in \cref{13}. The reconstructed profile demonstrates good agreement with the original $\sin(x)$ function.

\begin{figure}[htbp]
    \centering
    \begin{subfigure}{0.325\textwidth}
        \includegraphics[width=\linewidth]{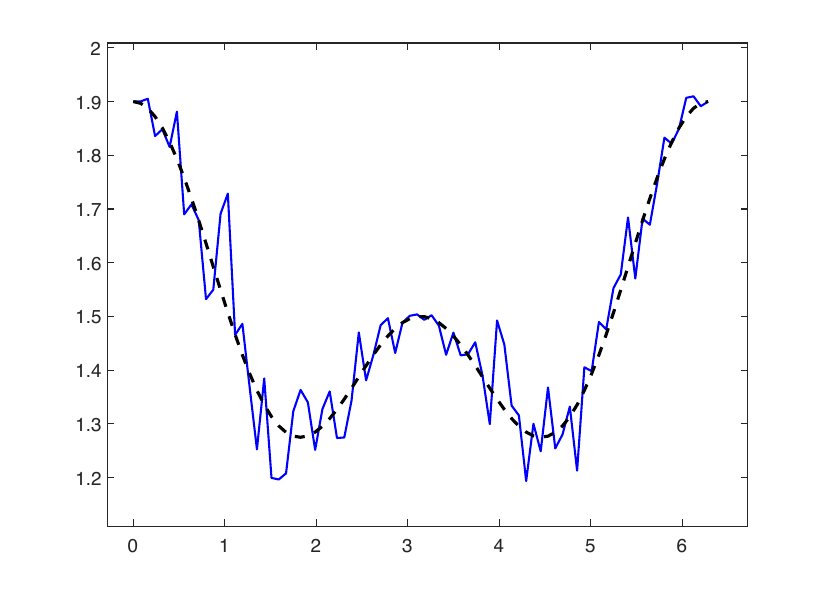}
        \caption{\it{\small   $f^{10}_{lip}(x)$ with $N_0=80$}}
        \label{180}
    \end{subfigure}
    \begin{subfigure}{0.325\textwidth}
        \includegraphics[width=\linewidth]{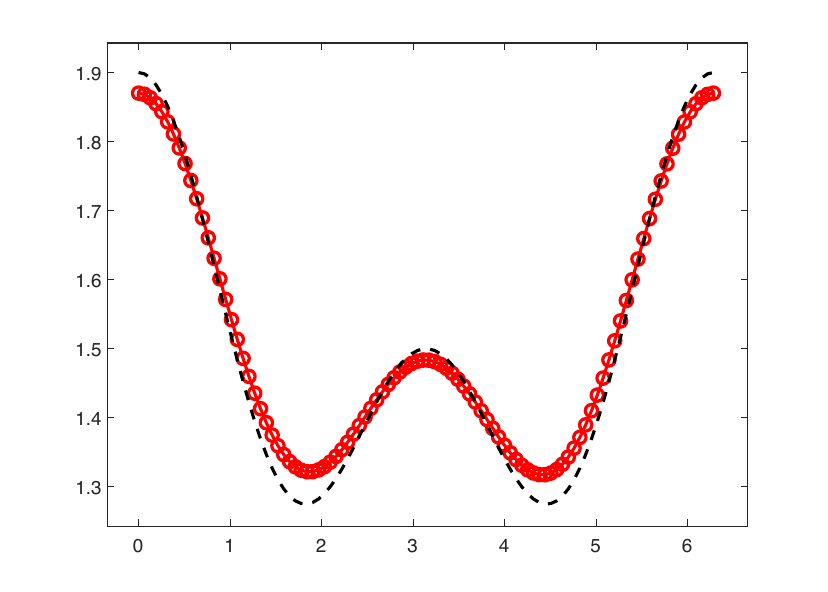}
        \caption{\it{\small   $\bar{f}(x)$ with $N_0=80$}}
        \label{180ginv}
    \end{subfigure}
       \begin{subfigure}{0.325\textwidth}
        \includegraphics[width=\linewidth]{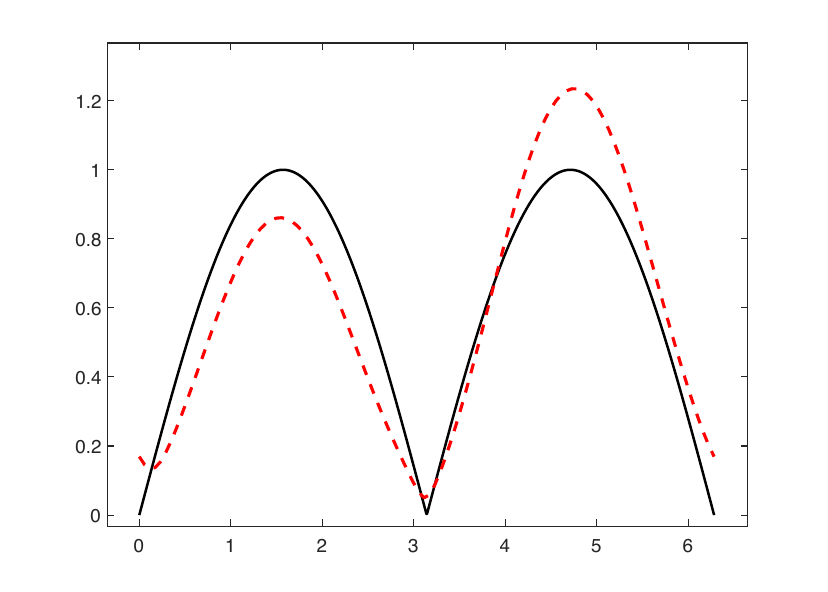}
        \caption{\it{\small  $|h_2(x)|$ with $N_0=80$}}
        \label{180inv}
    \end{subfigure}

     \begin{subfigure}{0.325\textwidth}
        \includegraphics[width=\linewidth]{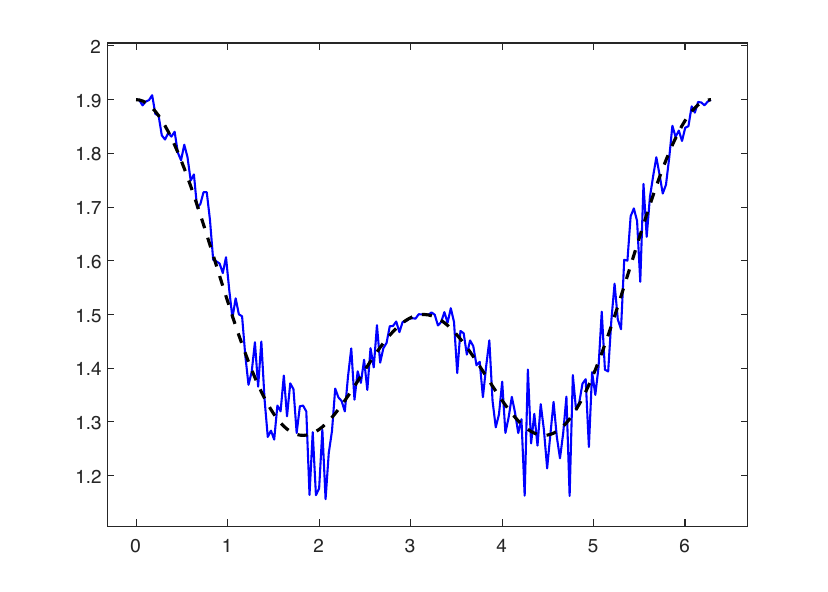}
        \caption{\it{\small $f^{10}_{lip}(x)$ with $N_0=160$ } }
        \label{1180}
    \end{subfigure}
\begin{subfigure}{0.325\textwidth}
        \includegraphics[width=\linewidth]{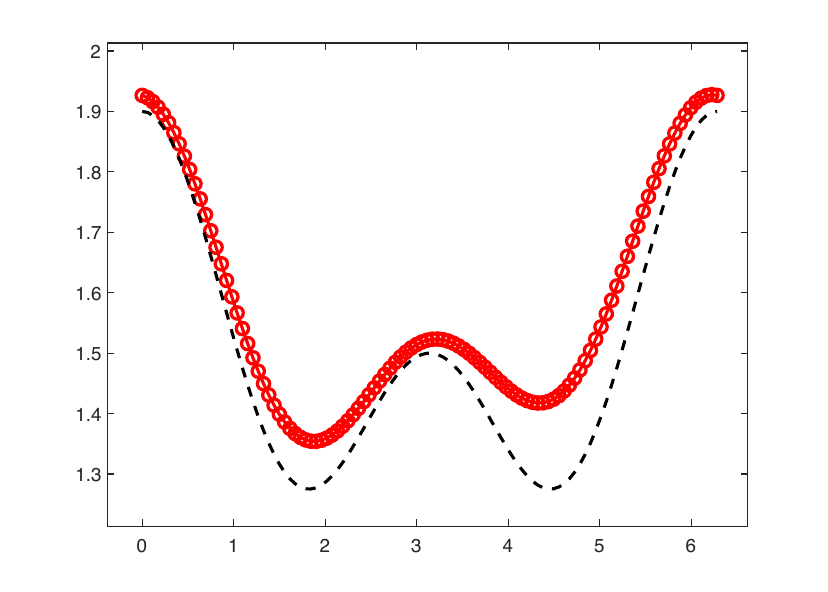}
        \caption{\it{\small $\bar{f}(x)$ with $N_0=160$ } }
        \label{1180ginv}
    \end{subfigure}
 \begin{subfigure}{0.325\textwidth}
        \includegraphics[width=\linewidth]{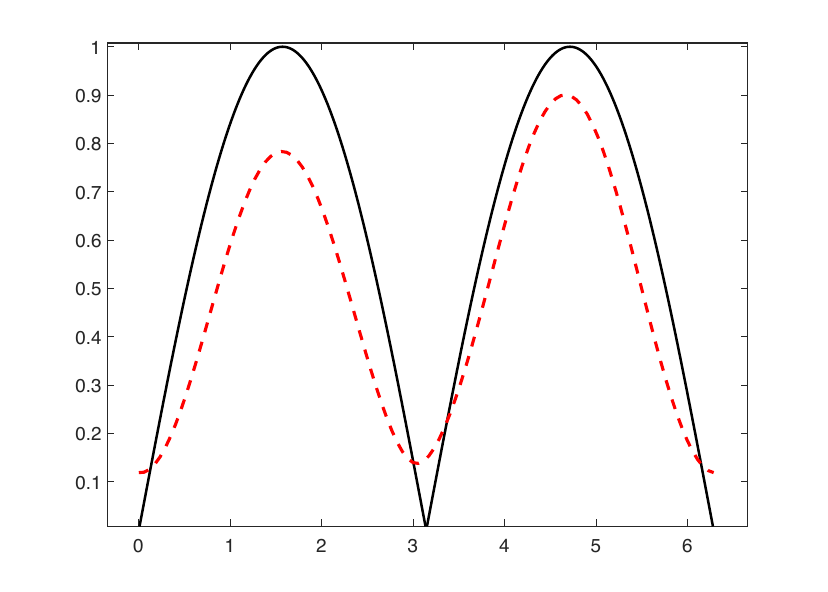}
        \caption{\it{\small  $|h_2(x)|$ with $N_0=160$ } }
        \label{1180inv}
    \end{subfigure}
    \caption{\it{\small The realization and the inversion at two different spatial resolutions for \Cref{ex2}. The top row corresponds to a coarse grid $N_0=80$ while the bottom row uses a finer grid $N_0=160$. First column:  the 10-th realization of random surface at two different spatial resolutions (Solid line) and  original $g_2(x)$ (Dashed line);  Second column: the reconstructed mean profile $\bar{f}(x)$ (Circular markers) compared with the ground truth $g_2(x)$(Dashed line).  Third column: the reconstructed variance intensity $|h_2(x)|$ (Dashed line)  compared with the true function $|sin(x)|$ (Solid line).}}
    \label{ex1inv}
    \end{figure}

The inversion performance is compared in \cref{ex1inv} for other two spatial resolutions ($\Delta x=80,160$), showing both the raw realizations and reconstructed $g_2(x)$, $|h_2(x)|$. These results confirm the stable reconstruction of $h_2(x)$ within acceptable error bounds for the examined $\Delta x$ range, validating the scheme's discretization invariance. 

\begin{figure}[htbp]
    \centering
    \begin{subfigure}{0.325\textwidth}
        \includegraphics[width=\linewidth]{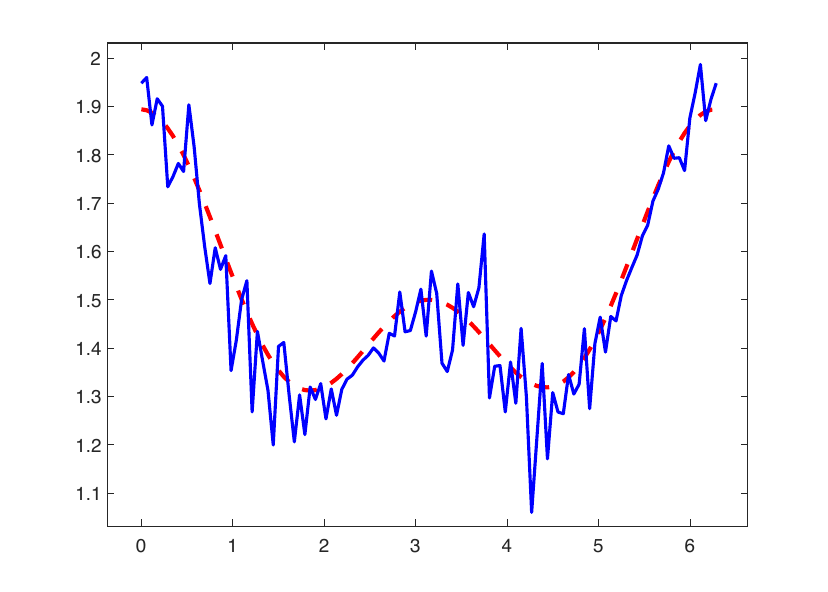}
        \caption{\it{\small Reconstruction $f_{10}(x)$}}
        \label{3re0}
    \end{subfigure}
     \begin{subfigure}{0.325\textwidth}
        \includegraphics[width=\linewidth]{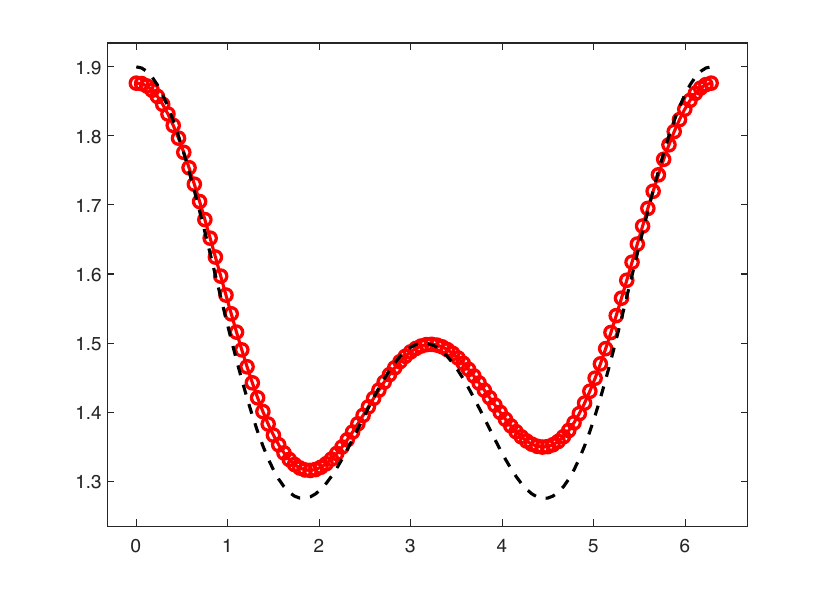}
        \caption{\it{\small Reconstructed $\bar{f}(x)$ } }
        \label{2ginnv}
    \end{subfigure}
     \begin{subfigure}{0.325\textwidth}
  \includegraphics[width=\linewidth]{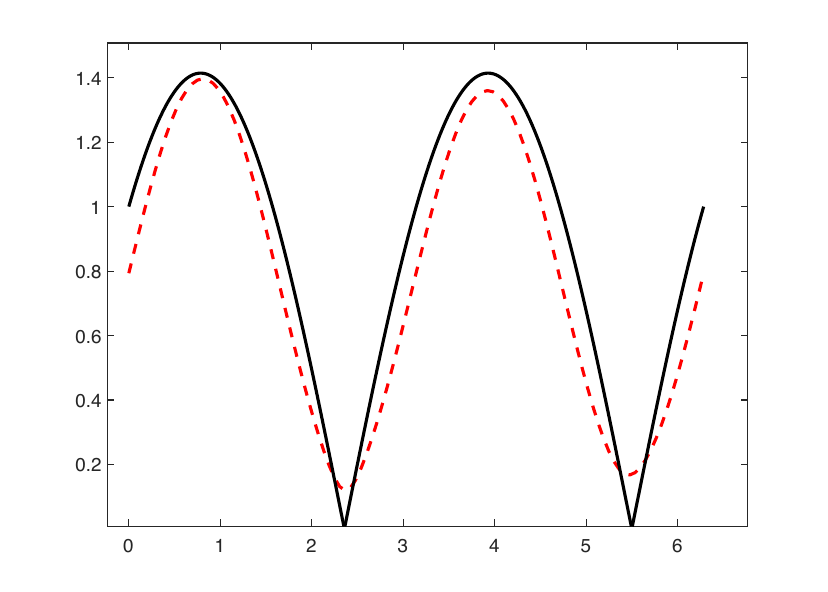}
  \caption{Reconstruted $h_3(x)$}
        \label{hsincos}
    \end{subfigure}
    \caption{\it{\small Inversion of $g_3(x)$ and $h_3(x)$. (a): Solid line: 10th surface $f^{10}_{lip}(x)$. Dashed line: the inversion result $f_{10}(x)$ obtained at the end of the tenth algorithm execution. (b): Dashed line: the original image of $g_3(x)$. Circular marked line line: the reconstructed $\bar{f}(x)$. (c): The inversion of $|h——3(x)|$ in \cref{ex3}. Solid line: the original image of $|\sin(x)|$. Dashed line: the  image of reconstructing $|h_3(x)|$.}}
    \label{ex33}  
    \end{figure}

\begin{example}{\label{ex3}}
     Reconstruct the random surface where
\[g_3(x)=1.5+0.2\cos(x)+0.2\cos(2x),\quad h_3(x)=\sin(x)+\cos(x),\] inside the interval $[0,2\pi]$.
\end{example}
Since $g_3(x)$ coincides exactly with $g_2(x)$ in \cref{ex2}, we adopt the representative resolution $N_0=110$ for demonstration purposes. The acquired dataset yields the inversion results displayed in \cref{ex33}.  In this particular case, the inversion outcomes associated with $g_3(x)$ and $h_3(x)$ exhibit a high degree of accuracy.

\begin{figure}[htbp]
    \centering
    \begin{subfigure}{0.325\textwidth}
        \includegraphics[width=\linewidth]{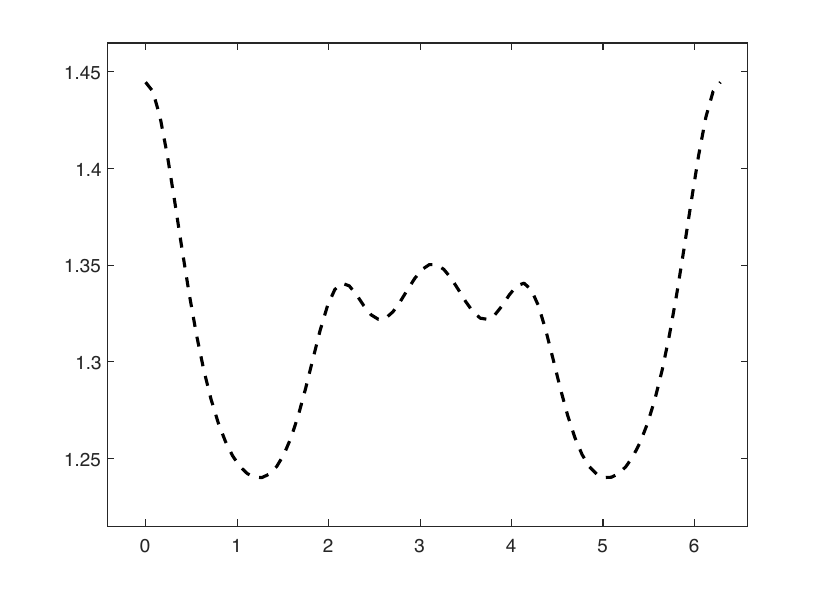}
        \caption{\it{\small Noise-free surface $g_4(x)$ }}
        \label{f45}
    \end{subfigure}
     \begin{subfigure}{0.325\textwidth}
        \includegraphics[width=\linewidth]{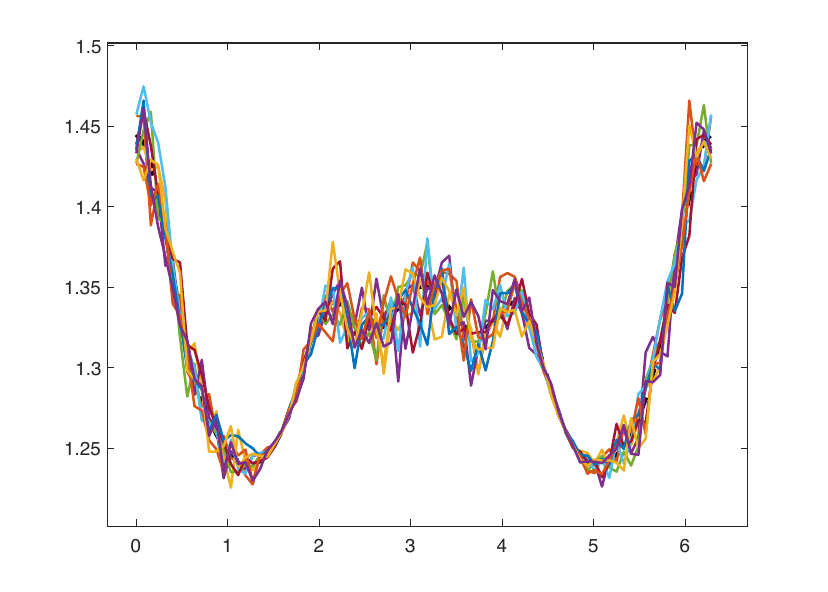}
        \caption{\it{\small Different realizations }}
        \label{flip45}
    \end{subfigure}
    \begin{subfigure}{0.325\textwidth}
        \includegraphics[width=\linewidth]{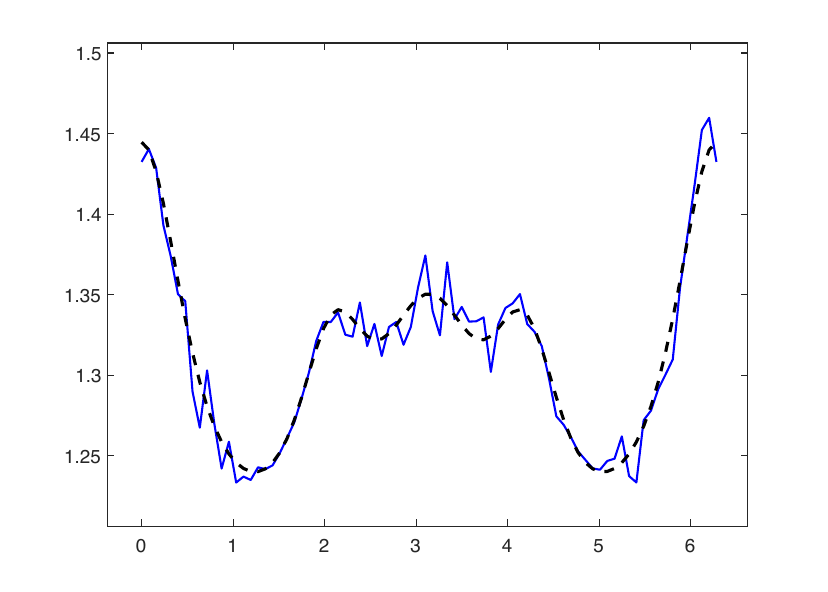}
        \caption{\it{\small The $10^{th}$ surface  $f^{10}_{lip}(x)$}}
        \label{f4510}
    \end{subfigure}
     \caption{\it{\small Numerical simulation of \cref{ex4} at spacial resolution $N_0=80$. (a): The image of noise-free function $g_4(x)$. (b): The realizations of first 10 samples $\{f^m_{lip}(x)\}_{m=1}^{10}$. (c): The $10$-th concrete realization $f^{10}_{lip}(x)$ of random surface.}}
    \label{exxx4}
\end{figure}

\begin{figure}[htbp]
    \centering
    \begin{subfigure}{0.325\textwidth}
        \includegraphics[width=\linewidth]{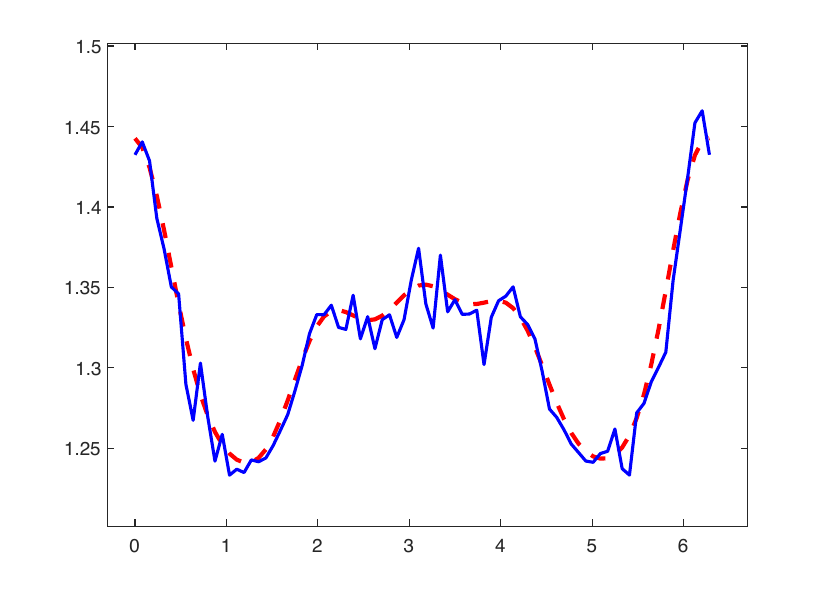}
        \caption{\it{\small Reconstruction $f_{10}(x)$}}
        \label{ffrei}
    \end{subfigure}
     \begin{subfigure}{0.325\textwidth}
        \includegraphics[width=\linewidth]{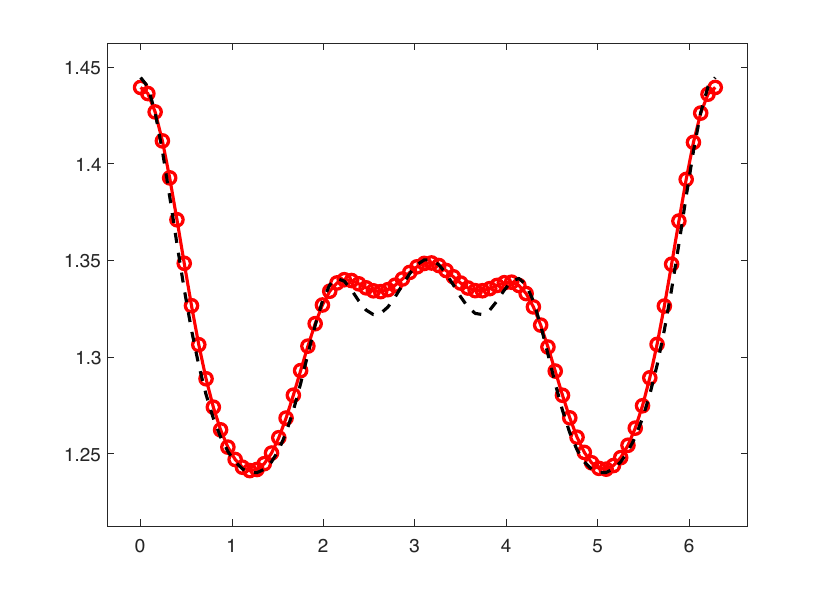}
        \caption{\it{\small Reconstructed $\bar{f}(x)$ } }
        \label{4ginv}
    \end{subfigure}
     \begin{subfigure}{0.325\textwidth}
  \includegraphics[width=\linewidth]{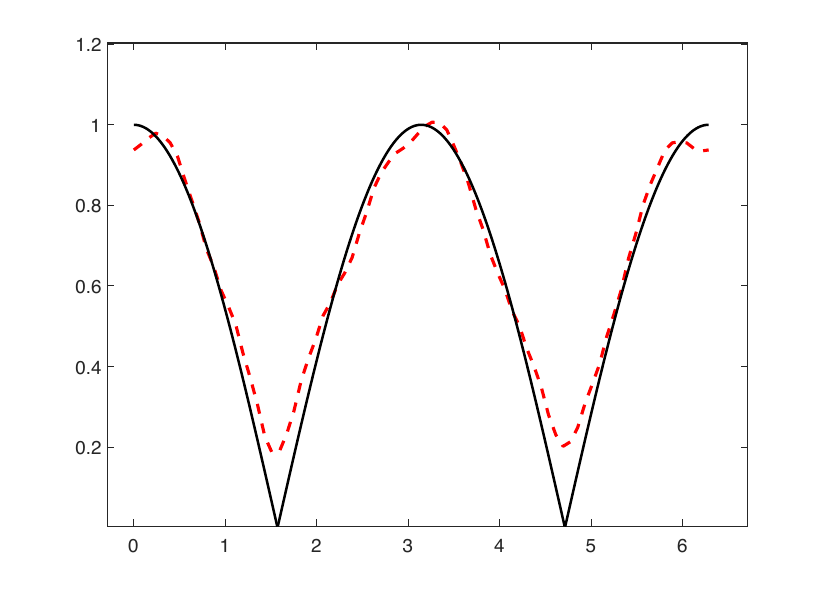}
  \caption{Reconstructed $|h_4(x)|$}
        \label{4hinv}
    \end{subfigure}
    \caption{\it{\small Inversion of $g_4(x)$ and $h_4(x)$. (a): Solid line: 10th surface $f^{10}_{lip}(x)$. Dashed line: the inversion result $f_{10}(x)$ obtained at the end of the tenth algorithm execution. (b): Dashed line: the original image of $g_4(x)$. Circular markers: the reconstructing $\bar{f}(x)$ . (c):The inversion of $|h_4(x)|$ in \cref{ex4}. Solid  line: the original image of $|\cos(x)|$. Dashed line: the  image of reconstructing $|h_4(x)|$.}}  
    \label{4exx}
    \end{figure}

\begin{example}{\label{ex4}}
    Reconstruct the random surface where
    \[g_4(x)=1.2+0.05e^{\cos(2x)}+0.04e^{\cos(3x)},\quad h_4(x)=\cos(x),\] inside the interval $[0,2\pi]$.
\end{example}
$g_4(x)$ in this case features a more complex functional form compared to the preceding examples. Consistent with the preceding analysis pipeline, we initiate this example presentation by visualizing the intuitive patterns of $g_4(x)$ and $f^{10}_{lip}(x)$ as shown in \cref{exxx4}.  The inversion of $g_4(x)$ achieves optimal accuracy at wavenumber $\kappa=6$. Set $k_{max}=6$.  With statistically sufficient measurements, the inversion outcomes illustrated in \cref{4exx} can be achieved. 
To enable direct comparison with our inversion results, $g_4(x)$ is represented by its Fourier series approximation, where the spectral coefficients  are compared against the reconstructed coefficients in \cref{tabel03}. Additional tests with $N_0=60$ and $N_0=110$ both yield accurate $|h_4(x)|$ inversion, though omitted here for brevity. 

\begin{table}[htbp]
  \centering
  \begin{tabular}{lccccccc}
    \toprule
    & $c_{m,0}$ & $c_{m,1}$ & $c_{m,2}$ & $c_{m,3}$ & $c_{m,4}$  & $c_{m,5}$  & $c_{m,6}$  \\
    \midrule
    Target\ \big($g_4(x)$\big)        & 1.3139   & 0 & 0 & 0.0565 & 0    & 0.0452 & 0     \\
    Initial   & 0   & 0 & 0 & 0 & 0    & 0 & 0       \\
    Computed\ \big($\bar{f}(x)$\big) & 1.3186   & 0.0008 & 0 & 0.0574 & -0.0002    & 0.0446 & 0    \\
    \bottomrule
     \toprule
   & $c_{m,7}$ & $c_{m,8}$ & $c_{m,9}$ & $c_{m,10}$ & $c_{m,11}$ & $c_{m,12}$ &   \\
    \midrule
    Target\ \big($g_4(x)$\big)    & 0.0136  & 0 & 0 & 0 & 0.0131   & 0 &     \\
    Initial  & 0   & 0 & 0 & 0 & 0    & 0 &       \\
    Computed\ \big($\bar{f}(x)$\big) & 0.0114   & 0.0004 & 0 & 0 & 0.0068    & 0.0002 &   \\
    \bottomrule
  \end{tabular}
   \caption{\it{\small Accuracy of reconstructed Fourier coefficients for $g_4(x)$ versus reconstructing $\bar{f}(x)=c_{m,0} + \sum_{p=1}^{2} [ c_{m,2p-1} \cos(px) + c_{m,2p} \sin(px) ]$ of \cref{ex4}. Target $g_4(x)$: Theoretical Fourier coefficients of the ground truth function. Initial: Initial guess for coefficients. Computed $\bar{f}(x)$: Reconstructed mean estimate  with $k_{max}=6$.}} 
  \label{tabel03}
\end{table}

\begin{example}{\label{ex5}}
     Reconstruct the random surface where
      \[g_5(x)=1.2+0.05e^{\cos(2x)}+0.04e^{\cos(3x)},\quad h_5(x)=\cos(x)+\cos(2x),\] inside the interval $[0,2\pi]$.
\end{example}

\begin{figure}[htbp]
    \centering
    \begin{subfigure}{0.325\textwidth}
        \includegraphics[width=\linewidth]{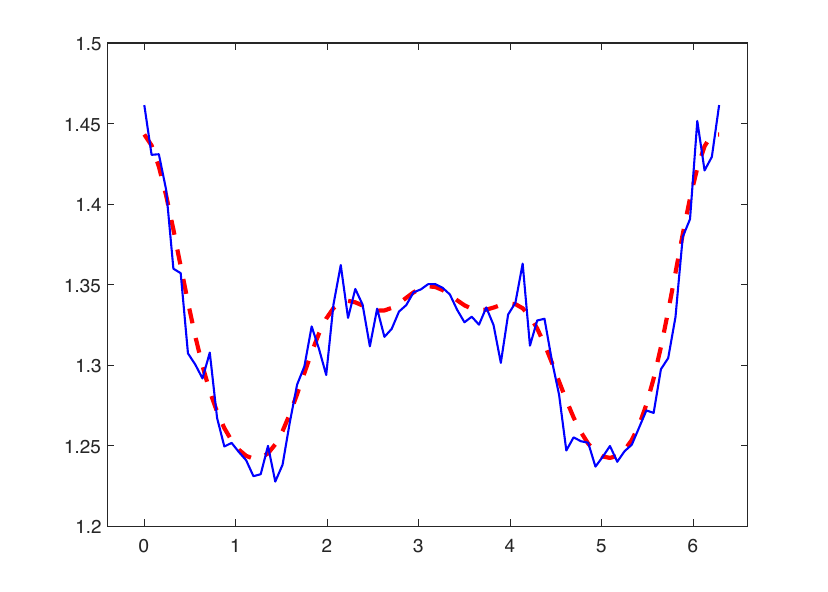}
        \caption{\it{\small Reconstruction $f_{10}(x)$}}
        \label{5ginv}
    \end{subfigure}
     \begin{subfigure}{0.325\textwidth}
        \includegraphics[width=\linewidth]{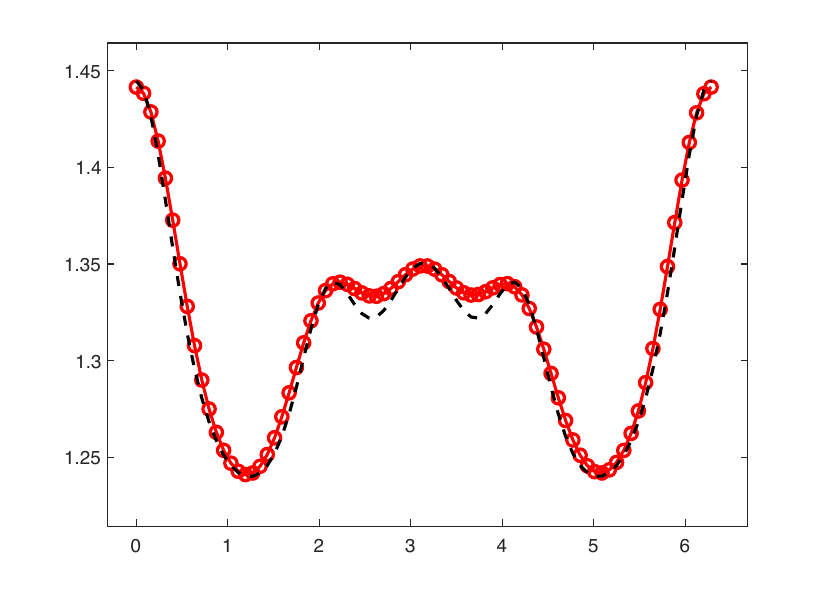}
        \caption{\it{\small Reconstructed $\bar{f}(x)$ } }
        \label{5ginnv}
    \end{subfigure}
     \begin{subfigure}{0.325\textwidth}
  \includegraphics[width=\linewidth]{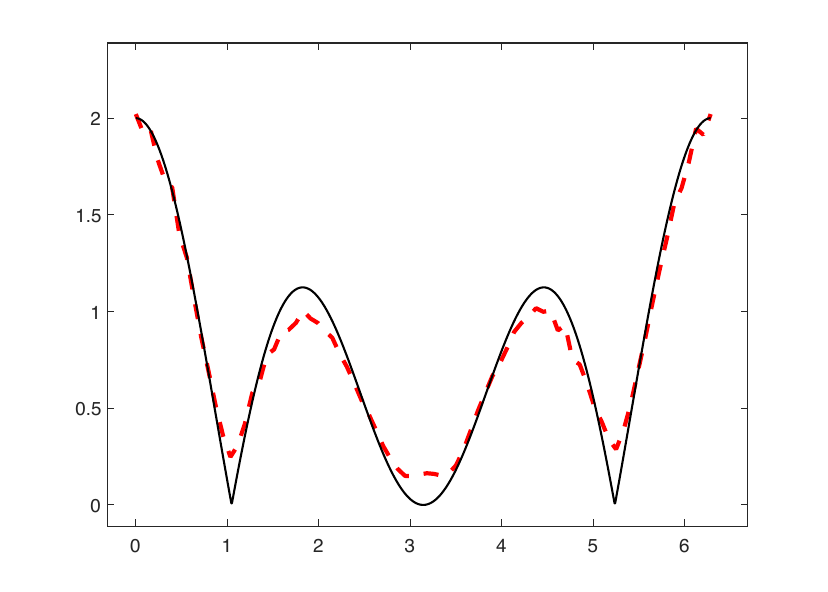}
  \caption{Reconstructed $|h_5(x)|$}
        \label{5hinnv}
    \end{subfigure}
    \caption{\it{\small Inversion of $g_5(x)$ and $h_5(x)$. (a): Solid line: 10th surface $f^{10}_{lip}(x)$. Dashed line: the inversion result $f_{10}(x)$ obtained at the end of the tenth algorithm execution. (b): Dashed line: the original image of $g_5(x)$. Circular markers: the reconstructing $\bar{f}(x)$. (c):The inversion of $|h_5(x)|$ in \cref{ex5}. Solid  line: the original image of $|\cos(x)+\cos(2x)|$. Dashed  line: the  image of reconstructing $|h_5(x)|$.}}  
    \label{5exx}
    \end{figure}

The final instance demonstrates the inversion performance for $h_5(x)=\cos (x) + \cos (2x)$. Visual inspection of the results in \Cref{5exx} confirms close agreement between the reconstructed and reference profiles. 
\section{Observation and discussion}{\label{sec:obs}}

\subsection{The analysis of inversion of $|h(x)|$}{\label{hxana}}

Based on the numerical examples presented in \Cref{sec:experiments},  we would like to analyze and discuss  several key observations   regarding the reconstruction performance.  We shall begin by identifying the critical factors governing the precision of the $|h(x)|$ estimation. Numerical results demonstrate that  the reconstruction accuracy of $|h(x)|$ exhibits strong dependence on the precision achieved in recovering $g(x)$. For  \cref{ex2} and \cref{ex3}, our results demonstrate that peak reconstruction accuracy for $g(x)$ occurs at $\kappa=2$ (see \cref{556}), which correspondingly yields an optimal $|h(x)|$ reconstruction as shown in \cref{13}. In the case where $\kappa\neq 2$, the deviation of the $\bar{f}(x)$ from the ground truth $g(x)$ will lead to  significant errors of inversion of $|h(x)|$. This dependent relation is further exemplified in \cref{ex4}: similarly,  optimal reconstruction of $g$ is uniquely achieved at $\kappa=6$ in \cref{ex4},  which in turn produces the most accurate $|h(x)|$ estimation. When $\kappa \neq 6$, the current methodology fails to maintain satisfactory reconstruction accuracy for $|h(x)|$, highlighting the algorithm's sensitivity to the accuracy of the inversion of $g(x)$.

Another crucial observation concerns the strong dependence of $|h(x)|$ reconstruction accuracy on the intensity of stochastic perturbations. As noted in \cref{911}, the amplitude of random fluctuations must remain moderate relative to the baseline profile. During our experimental investigation, it was found that this constraint not only inherently limits the achievable precision in reconstructing $|h(x)|$ but also restricts the method's capacity to recover more complex forms of $h(x)$.  Excessive perturbation amplitudes may enable the handling of intricate $h(x)$, but  it will also  severely distort the grating's fundamental topography $g(x)$,  resulting in significant errors in $g(x)$ reconstruction. Moreover, overly intense disturbances will produce physically implausible surface features in real-world scenarios. Consequently, a delicate trade-off emerges: precise reconstruction of $|h(x)|$ requires accurate inversion of $g(x)$, yet efforts to enhance $g(x)$ precision may compromise $|h(x)|$ accuracy. By identifying the optimal quantitative balance between these factors, the desired reconstruction outcome can be achieved. If we have some prior knowledge about the general shape of $h(x)$, we can identify the nodes that have been affected by this interference based on the trend of sign changes in the majority of samples and make corrections accordingly, thereby obtaining an accurate inversion of $h(x)$.

Our numerical experiments in \Cref{sec:experiments} demonstrate that the inversion is highly effective when $h(x)$ is a linear combination of low-frequency trigonometric functions. Attempts to reconstruct more complex, non-trigonometric functions, however, were unsatisfactory under the same conditions. We speculate that this performance disparity is due to the efficient phase matching of low-frequency trigonometric functions with propagating wave modes dictated by the Rayleigh expansion. In contrast, high-frequency information—corresponding to the fine spatial details in both high-frequency trigonometric terms and more complex, non-trigonometric functions—is difficult to capture. This is because these features predominantly excite evanescent waves, which decay exponentially before reaching the detector.  A more rigorous mathematical analysis establishing these connections is an interesting direction for future research.

\subsection{Possibility of reconstructing $h(x)$}{\label{hxinv}}
This section we are going to discuss the possibility of precisely reconstructing $h(x)$. At the final algorithmic stage, the immediate outcome from \cref{MCCh}  is $h^2(x)$ —— thus the core challenge of recovering $h(x)$ lies in how to resolve the sign ambiguity between ``$h(x)$" and ``$-h(x)$" from $|h(x)|$. Notably,  the presence of noise $\xi_i$ may introduce additional complexity,  potentially causing sign alternation where both  $h(x)$ and $-h(x)$ appear inconsistently in the reconstructed output. Taking \cref{ex1}  as an illustrative case. To recover $\cos(x)$ from $|\cos(x)|$,  we will typically encounter four distinct curve patterns in the end: $\cos(x)$, $-\cos(x)$, $M$-shaped and $W$-shaped curves. Fortunately, from a probabilistic perspective, the chance of $M$-shaped and $W$-shaped curves occurring is considerably less than that of $\cos(x)$ and $-\cos(x)$ patterns, particularly as the number of calibration samples grows (See experiment results in \cref{ingv}). Based on these observations,  we propose a methodological framework for accurate inversion of $h(x)$: While computing $c_{ij}$ in \cref{hin}, record the sign of $ f_m(x_i)-\mathbb{E}f(x_i)$ at each nodal point $x_i$ simultaneously. This procedure will ultimately yield  a sign vector $\text{sgn}(h)=\{+,-\}^{N_0}$. At the final step of  inversion, the sign of $h(x_i)$ will be recovered using the recorded information $\text{sgn}(h)$. 
 In the illustrated \cref{ex1}, applying our previously discussed approach and plotting the results, we obtain an image depicted in \cref{hinv}:
\begin{figure}[htbp]
    \centering
    \begin{subfigure}{0.5\textwidth}
        \includegraphics[width=\linewidth]{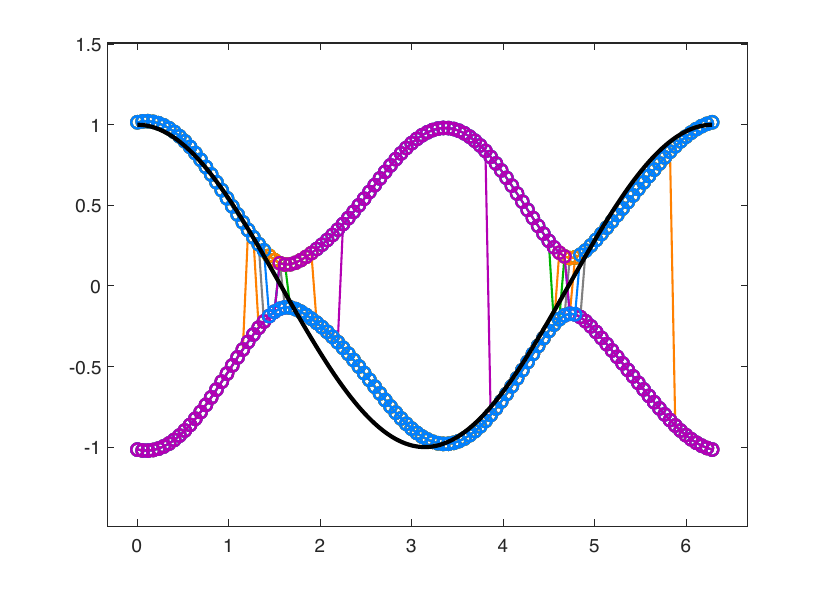}
        \caption{\it{\small Reconstructed $h(x)$ with sign information}}
        \label{hinv}
    \end{subfigure}
     \begin{subfigure}{0.5\textwidth}
        \includegraphics[width=\linewidth]{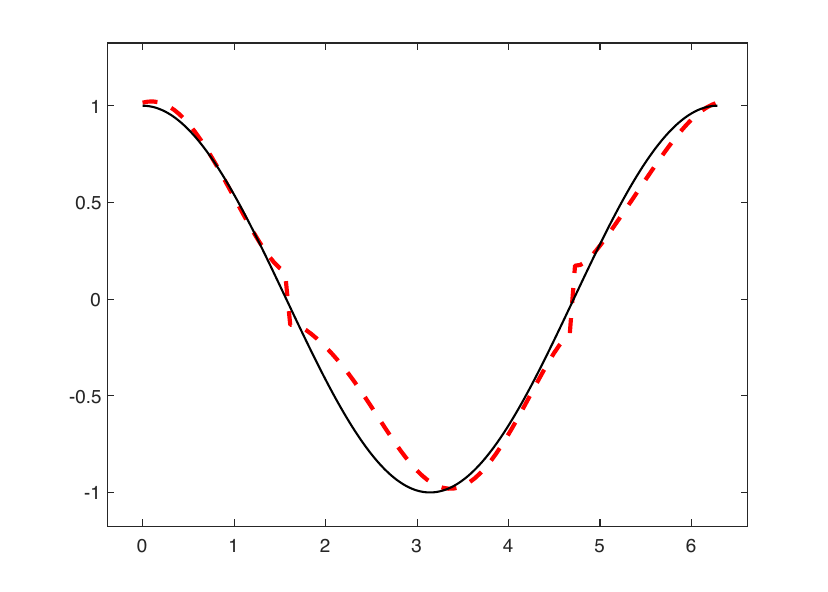}
        \caption{\it{\small Ultimate reconstructed $h(x)$ } }
        \label{hinvv}
    \end{subfigure}
    \caption{\it{\small Accurate inversion of $h(x)=\cos(x)$. (a): Purple line: `$-\cos(x)$' pattern of inversion  with 90 samples. Blue line: `$\cos(x)$' pattern of inversion with 98 samples. (b): Solid line: the original image of $\cos(x)$. Dashed line: the reconstructed $h(x)$ obtained using statistical analysis and prior information.}}
    \label{ingv}
    \end{figure}
We observed that among the first 200 samples, 90 samples (45\%) exhibited  a ``$-\cos(x)$'' pattern (Purple line); 98 samples (49\%) displayed the correct ``$\cos(x)$'' pattern (Blue line); and 12 samples (6\%) formed patterns in the shape of ``$W$" and ``$M$". Despite minor sign fluctuations at certain nodes (visible as vertical line segments in \cref{hinv}), the substantial number of supporting samples suggests that we can reliably recover the two dominant patterns: ``$\cos(x)$'' and ``$-\cos(x)$''. Furthermore, when incorporating some other prior knowledge about $h(x)$ such as ``$h(0)>0$", we can accurately determine the correct inversion of $h(x)$, as demonstrated in \cref{hinvv}.

Since the function $h(x)$ we consider in this paper belongs to a class of simple, well-behaved functions, their intrinsic boundedness and symmetry properties create the possibility for us to utilize some  statistical approaches for precise inversion of $h(x)$.  As the  complexity of shape structure to $h(x)$ escalates, a proliferation of patterns that deviate from the original $h(x)$ and $-h(x)$ configurations may emerge successively. This phenomenon, in turn, significantly amplifies the  difficulty  and challenge associated with the inversion process of $h(x)$. Our proposed method demonstrates feasibility  for simple $h(x)$ functions, it remains an initial step toward handling more general cases, leaving room for further methodological developments. Extensions to more complex cases await future research.

\section{Conclusions and Future work}{\label{sec:conclusions}}

In this paper,  we have investigated the inverse scattering problem for  a time-harmonic plane wave incident on a perfectly reflecting random periodic structure. We introduced a novel modeling approach where each realization of the random surface is a Lipschitz-continuous diffraction grating, constructed via the discretization of a stochastic process and linear interpolation.  Mathematically, the surface profile is expressed as the sum of a baseline profile and a linear combination of weighted local `tent' basis functions. To reconstruct the core statistics of this random grating, we propose the Recursive Parametric Smoothing Strategy (RPSS). By integrating this strategy with Monte Carlo sampling and a wavenumber continuation approach, our reconstruction scheme MCCh demonstrates effectiveness under certain elementary benchmark configurations. 

Note that there are several promising avenues for future exploration. First, a rigorous theoretical analysis of the proposed algorithm should be conducted. Second, an extension for more realistic grating structures which often involves non-uniform node distributions should be investigated. Finally, advancing the algorithm’s capability to achieve accurate reconstruction of more complex statistical functions, particularly under conditions of limited sample support, would be a significant and challenging direction for future research. Looking ahead, we are going to explore more robust approaches with reduced sample requirements in subsequent research, with the goal of reconstructing key statistical quantities for more complex random surface configurations.

\section*{Acknowledgments}
This work was partially supported by National Key R\&D Program of China (no. 2024YFA1016100) and National Natural Science Foundation of China (no. U21A20425, no. 12201404). The computations in this paper were run on the Siyuan-1 cluster supported by the Center for High Performance Computing at Shanghai Jiao Tong University. We thank Prof. Gang Bao for helpful advice on the manuscript.

\appendix

\bibliographystyle{siamplain}
\bibliography{accept}
\end{document}


\maketitle

\section{A detailed example}

Here we include some equations and theorem-like environments to show
how these are labeled in a supplement and can be referenced from the
main text.
Consider the following equation:
\begin{equation}
  \label{eq:suppa}
  a^2 + b^2 = c^2.
\end{equation}
You can also reference equations such as \cref{eq:matrices,eq:bb} 
from the main article in this supplement.

\lipsum[100-101]

\begin{theorem}
  An example theorem.
\end{theorem}

\lipsum[102]
 
\begin{lemma}
  An example lemma.
\end{lemma}

\lipsum[103-105]

Here is an example citation: \cite{KoMa14}.

\section[Proof of Thm]{Proof of \cref{thm:bigthm}}
\label{sec:proof}
\lipsum[106-112]

\section{Additional experimental results}
\Cref{tab:foo} shows additional
supporting evidence. 

\begin{table}[htbp]
{\footnotesize
  \caption{Example table}  \label{tab:foo}
\begin{center}
  \begin{tabular}{|c|c|c|} \hline
   Species & \bf Mean & \bf Std.~Dev. \\ \hline
    1 & 3.4 & 1.2 \\
    2 & 5.4 & 0.6 \\ \hline
  \end{tabular}
\end{center}
}
\end{table}

\bibliographystyle{siamplain}
\bibliography{references}